\documentclass[english]{article}

\usepackage{natbib}
 \bibpunct[, ]{(}{)}{,}{a}{}{,}%

\usepackage{algorithmicx,algorithm}
\usepackage[noend]{algpseudocode}
\usepackage{graphicx}
\usepackage{amsmath}
\usepackage{multirow}
\usepackage{subfig}
\usepackage{amsmath}
\usepackage{amssymb}
\usepackage{setspace}
\usepackage[T1]{fontenc}
\usepackage[latin9]{inputenc}
\usepackage{geometry}
\usepackage{babel}
\usepackage{mathrsfs}
\usepackage{mathtools}
\usepackage{amsmath}
\usepackage{amsthm}
\usepackage{amssymb}
\usepackage{authblk}
\usepackage{breakurl}
\usepackage{algorithm}
\usepackage{algpseudocode}
\usepackage{tikz}
\usepackage{enumerate} 
\usepackage{graphicx}  
\usepackage{wrapfig} 
\usepackage{caption} 
\usepackage{hyperref}

\usepackage{natbib}
\bibpunct[, ]{(}{)}{,}{a}{}{,}%

\makeatletter
\theoremstyle{plain}
\newtheorem{thm}{\protect\theoremname}[section]
\theoremstyle{definition}
\newtheorem{defn}[thm]{\protect\definitionname}
\theoremstyle{plain}
\newtheorem{assumption}[thm]{\protect\assumptionname}
\theoremstyle{plain}

\theoremstyle{remark}

\theoremstyle{plain}

\theoremstyle{plain}
\newtheorem{lem}[thm]{\protect\lemmaname}

\makeatother

\providecommand{\assumptionname}{Assumption}
\providecommand{\corollaryname}{Corollary}
\providecommand{\definitionname}{Definition}
\providecommand{\lemmaname}{Lemma}
\providecommand{\propositionname}{Proposition}
\providecommand{\remarkname}{Remark}
\providecommand{\theoremname}{Theorem}

\begin{document}





\title{An Online Non-Stationary Simulation Optimization Approach Based on Regime Switching}


\author[1]{Jianglin Xia}
\author[2]{Haowei Wang}
\author[1,3]{Songhao Wang}
\author[2]{Szu Hui Ng}
\affil[1]{College of Business, Southern University of Science and Technology, Shenzhen, CHINA}
\affil[2]{Department of Industrial Systems Engineering and Management, National University of Singapore, Singapore, SINGAPORE}
\affil[3]{Corresponding author: wangsh2021@sustech.edu.cn}
\date{}

\maketitle

\begin{abstract}
Dynamic and evolving operational and economic environments present significant challenges for decision-making. We explore a simulation optimization problem characterized by non-stationary input distributions with regime-switching dynamics across sequential decision stages. This problem encompasses both prediction uncertainty, arising from the regime-switching behavior of input distributions, and input uncertainty, resulting from parameter estimation for these distributions and their dynamics using finite data streams. To address these uncertainties, we develop a Bayesian framework that approximates the true objective function using a Markov Switching Model (MSM). We rigorously validate this approximation by establishing the consistency and asymptotic normality of the objective functions and optimal solutions. To solve the problem in an online fashion, we propose a metamodel-based algorithm that leverages simulation results from previous stages to enhance decision-making. Furthermore, we tackle scenarios with an unknown number of regimes through a Bayesian nonparametric method. Numerical experiments demonstrate that our algorithm achieves superior performance and robust adaptability.
\end{abstract} 





%

\section{Introduction}\label{sec:Intro}

The operational and economic environments of many real-world systems are dynamic and non-stationary, driven by gradual trends, regime shifts, or abrupt changes. This presents significant challenges for effective decision-making. For example, global supply chains can be affected by disturbances such as natural disasters, labor disputes, or infrastructure failures \citep{wilson2007impact}. Similarly, customer demand patterns fluctuate significantly due to large-scale pandemics, extreme weather events, or promotional activities \citep{keskin2022data, trapero2015identification}. Among various forms of non-stationarity, regime-switching is prominent and widely observed in real-world applications, where time-varying environment is driven by shifting regimes. A classic example is found in financial markets, where bull and bear market regimes-characterized by sustained positive or negative returns-necessitate distinct investment strategies. Additionally, the four-regime framework \citep{Bridgewater2012}, defined by combinations of rising or falling inflation and economic growth, illustrates how varying macroeconomic conditions uniquely influence the performance of diverse financial assets. Similarly, in cyclical category buying, regime shifts between high and low purchase tendencies enable more effective promotional strategies through precisely timed price adjustments \citep{park2011regime}.
These examples collectively highlight the need for decision-making approaches to navigate non-stationary and volatile environments, particularly in today's context of frequent political, economic, and environmental changes. 

To address this, the following stochastic optimization problem can be considered:
\begin{equation}
	\label{SO_problem}
	\min_{x\in \mathcal{X}} \; \mathbb{E}_{\xi\sim P^c}[y(x,\xi)],
\end{equation}
where $x$ is the decision variable from a compact space $\mathcal{X}$, and $\xi$ is the random environmental variable following distribution $P^c$. In the non-stationary settings above, $P^c$ evolves over time, such as the shifting demand distributions in inventory management and the regime-switching return distributions in finance.        
In this work, we consider problem \eqref{SO_problem} through \emph{\textbf{simulation optimization}} in \emph{\textbf{regime-switching}} contexts. Simulation optimization tackles problems with no closed-form and evaluates them using simulation models. This approach has emerged as a powerful tool when the real-world systems are complex or expensive to run \citep{hong2009brief}, including those in potential non-stationary environments. For example, in inventory control problems, the long-run expected cost per period is often analytically intractable, prompting the use of simulation to estimate and minimize the average cost \citep{fu1997techniques}. Similarly, in portfolio risk management, the absence of the closed-form portfolio loss distributions necessitates simulation to estimate risk measures and optimize the investment strategies \citep{hong2014monte}. 

\begin{figure}
	\centering
	\includegraphics[width=0.4\textwidth]{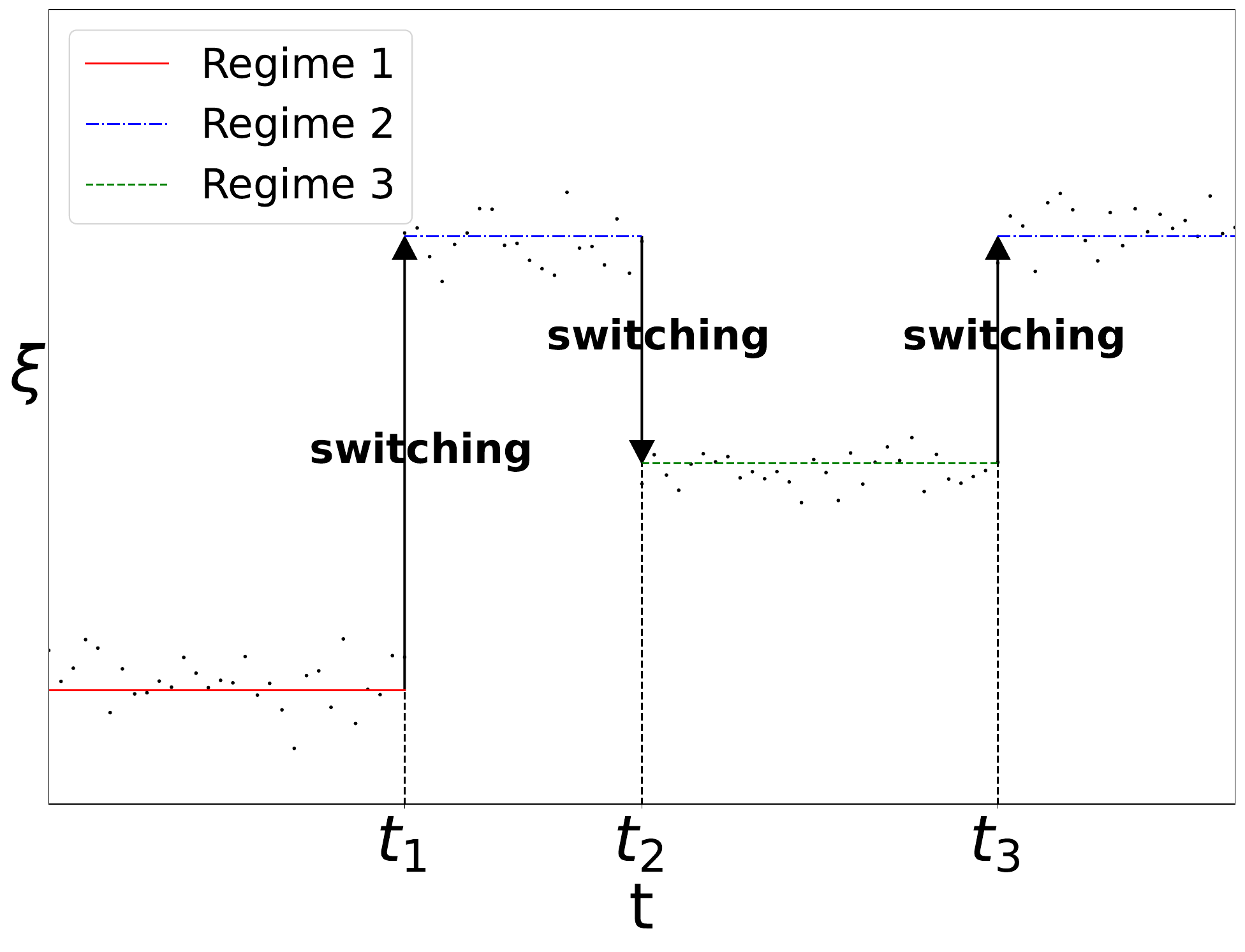}
	\caption{Illustration of regime switching process.}
	\label{fig:Regime}
\end{figure}

In a regime-switching process, the statistical properties of $P^c$, such as mean and volatility, vary across underlying regimes. Figure~\ref{fig:Regime} provides of an illustration, where $P^c$ within each regime is a unit-variance Gaussian but exhibits distinct means across regimes (indicated by horizontal lines).
Three regime switchings occur at time points $t_1\sim t_3$. Such regime changes are frequently observed in applications such as finance and inventory management, as noted earlier.
Throughout this work (except Section \ref{sec:Extension}), we assume the number of regimes is known a priori, a reasonable assumption supported by domain knowledge in many practical applications \citep{pun2023data,liu2022financial}. We further assume that both the distribution parameters and the switching points are unknown and should be inferred from historical data of $\xi$. When the historical data, such as the daily customer demand data during a selling season, are available and updated along multiple decision stages, we are expected to develop online simulation optimization strategies to leverage these data streams. Thus, the problem addressed in this work can be summarized as follows: within a decision period comprising multiple stages, observations of $\xi$ from prior stages are received at the start of each stage. The task is to predict the regime of the upcoming stage, estimate the unknown distribution $P^c$, and propose a decision through an efficient simulation optimization approach.  

In the simulation community, the random variable $\xi$ and its distribution $P^c$ are called the input variable and input distribution, respectively. Extensive studies have explored to approximate \eqref{SO_problem} with estimations of $P^c$, from historical data of $\xi$, and assess the approximation quality \citep{xie2014bayesian,zhou2015simulation, wang2020gaussian}. A key finding from these studies is that ignoring the input uncertainty, uncertainty from finite-data estimates of $P^c$, will severely compromise the solution quality (see a review in Section \ref{subsec:SO_IU}). In addition to input uncertainty, our problem introduces an additional layer of uncertainty: prediction uncertainty, arising from predicting the regime for the upcoming decision stage. Unlike input uncertainty that typically decreases as more input data become available, the prediction uncertainty here will persist due to abrupt regime switchings. These intertwined uncertainties pose substantial challenges for simulation optimization. 

To address these uncertainties, we develop a Regime-Switching Online Bayesian Simulation Optimization (RSOBSO) approach. To handle the prediction uncertainty, we employ a Markov Switching Model (MSM) to capture the regime-switching dynamics of non-stationary streaming input data. In MSM, a Markov chain is used to model the possible transitions among different regimes. It enables the construction of a predictive distribution for the regime in the next time stage. These predictive probabilities are then used to build a weighted objective function, reflecting the contribution of each possible future regime to the decision making. To account for the input uncertainty, we adopt Bayesian parameter estimation for the input distribution and the MSM dynamics, and approximate the objective function with its posterior. Our approximation extends the approaches of \cite{zhou2015simulation} and \citet{wu2018bayesian} by considering the joint uncertainty in both regime prediction and distribution parameter estimation.

The key contributions of this paper are summarized as follows. We formulate a Bayesian objective function as a principled approximation to the true objective, the one with known parameters for the input distribution and regime-switching dynamics, for the forthcoming decision stage. The approximation captures both prediction and input uncertainties. Our primary theoretical contributions lie in demonstrating the effectiveness of the proposed approximation through establishing the consistency and asymptotic normality of both the objective function value and optimal solutions. Specifically, as the number of time stages and input data size tend to infinity, the objective value and optimal solutions of the approximation converge to their counterparts under the true objective. In addition, our asymptotic normality results show that $\sqrt{t}$ ($t$ is the number of observations of $\xi$) times the difference between our estimated optimal values/solutions and the true ones follows a normal distribution asymptotically as the real world data size increases, revealing a rate of convergence at $O_p(1/\sqrt{t})$. Asymptotic normality results can be used for statistical inference (e.g., hypothesis testing and building confidence intervals) of the optimal solutions and objective.

Our computational contributions involve developing a metamodel-based algorithm to solve the proposed online problem. A key feature of the algorithm is the construction of a unified metamodel that jointly incorporates both decision variables and input parameters, enabling efficient reuse of simulation results from past decision stages to inform future decisions. Reusing simulation results is a critical issue for online simulation optimization, especially under the limited computing resource per stage \citep{wu2024data}. The results from previous stages, though run on different input models, carry valuable information for decision making in the current stage. With the joint model, we build a regime-aware Expected Improvement (EI) acquisition function to sequentially select the next design point and input distribution to run the experiment, aiming to efficiently guide the search for the optimal solution. To validate the algorithm both theoretically and empirically, we establish its theoretical convergence and demonstrate its strong empirical performance through numerical experiments. These include an inventory problem and a portfolio optimization task tested on the real-world regime and return data from the 2008-2009 Great Recession, a period extensively studied in financial research \citep{pun2023data, wang2025limited}.

The remainder of the paper is organized as follows. Section~\ref{sec:Literature} reviews the related works. Section~\ref{sec:Method} formulates the Bayesian approximated objective function. Section~\ref{sec:Theory} provides theoretical validations of the approximation with consistency and asymptotic normality results. The proposed RSOBSO algorithm is detailed in Section~\ref{sec:Algorithm} with convergence analysis. Section~\ref{sec:Experiment} presents numerical experiments. Section~\ref{sec:Extension} extends the approach to scenarios with unknown number of regimes through a Hierarchical Dirichlet Process. Finally, Section~\ref{sec:Conclusion} concludes the paper.

\section{Related Works}\label{sec:Literature}

We review relevant works in three aspects: simulation optimization under input uncertainty, simulation optimization with streaming input data, and regime-switching streaming data modeling. 

\subsection{Simulation Optimization Under Input Uncertainty}\label{subsec:SO_IU}
Input uncertainty stems from estimating input distributions with finite data and includes parametric uncertainty (estimating parameters within a known distribution family) and model uncertainty (related to the choice of distribution form) \citep{lam2016advanced,xie2014bayesian,barton2014quantifying}. Neglecting input uncertainty can lead to suboptimal designs and system risks. Consequently, researchers have increasingly focused on incorporating it into simulation optimization. One research stream focuses on discrete optimization, particularly ranking and selection (R\&S) \citep{wu2017ranking, xiao2018simulation, fan2020distributionally}. Another addresses continuous optimization. \citep{zhou2015simulation, wu2018bayesian} propose a Bayesian framework to hedge against the risk in parametric input uncertainty. To solve simulation optimization under this formulation, GP-based optimization algorithms have been adapted to handle input uncertainty effectively \citep{pearce2017bayesian, wang2020gaussian, cakmak2020bayesian}. Additionally, extensions have been made to nonparametric Bayesian models when the input distribution family is unknown \citep{wang2020nonparametric}.

\subsection{Simulation Optimization with Streaming Input Data}
In simulation modeling, it is typically assumed that a fixed batch of input data is available before decision-making. An alternative setting involves streaming input data, where new data arrive over time. Recent studies have explored simulation optimization in this setting. For example, \citet{song2019stochastic} 
apply a stochastic approximation framework to solve a sequence of evolving problems. For the R\&S problem with streaming input data, \citet{wu2024data} propose a fixed confidence method, while \citet{wang2025ranking} consider varying budgets for collecting multiple input data streams and develop the fixed budget formulation.
\citet{liu2024bayesian} address a decision-dependent problem, where the input distribution depends on the decisions, using an online Bayesian approach. 
These studies, however, have not tackled the challenge of input distribution exhibiting regime-switching non-stationarity driven by external factors.

\subsection{Regime-Switching Streaming Data Modeling}
Common non-stationary behaviors in streaming data have been summarized in \citet{iquebal2018change}, including varying moments, piecewise stationarity (regime-switching), and arbitrary variations. They also discuss modeling regime evolution as a discrete-time Markov chain. 
The classical MSM, also known as the Hidden Markov Model (HMM), is one such model. It has been widely applied in various fields due to its flexible nature and mathematical structure, including stock market, speech recognition, data processing, and bioinformatics \citep{mor2021systematic}. In the context of portfolio optimization, MSM has attracted considerable attention in recent years. However, most existing studies incorporate estimated MSM parameters directly into models without accounting for the associated estimation uncertainty \citep{costa2019risk, oprisor2020multi}. To address this issue, \citet{pun2023data} incorporate MSM into the distributionally robust optimization framework by constructing regime-switching ambiguity sets for the random returns, which, however, can be overly conservative and does not apply to simulation optimization problems. 

\section{Problem Formulation}\label{sec:Method}
In this section, we first illustrate the preliminaries of MSM in Section~\ref{sec:msm}, and then present the Bayesian formulation of the objective function and its equivalent form in Section~\ref{sec:obj}.   

\subsection{The Preliminaries of MSM}\label{sec:msm}
MSM $\{(S_t, \xi_t), t\geq1\}$ is a discrete-time random process formed by a Markov chain $\{S_t\}$ with finite regime space $S:=\{1,2,\cdots,\tilde{R}\}$ and fixed transition matrix $A:= (A_{i,j})_{i,j=1}^{\tilde{R}}$, and a sequence of random variables $\{\xi_t\}$. In a hidden MSM, $S_t$ is a latent variable indicating the regime of $\xi_t$. Each regime $s$ is associated with a unique emission distribution $\tilde{P}(\xi|\lambda_s)$ with parameter $\lambda_s$. Given $S_t$, $\xi_t$ follows $\tilde{P}(\xi|\lambda_{S_t})$ and is independent with other random variables. In Figure~\ref{fig:Regime}, the number of regimes is $\tilde{R}=3$, where each regime's emission distribution is a Gaussian with distinct means. The non-stationarity of $\xi$ arises from transitions among these three regimes. Such piecewise stationary behavior in data can be effectively captured and modeled using MSM. The random input parameter vector $\vartheta$ can thus be explicitly represented as $(\lambda_1,\cdots, \lambda_{\tilde{R}},A_{11}, \cdots,A_{\tilde{R}\tilde{R}})$ with a total dimension of $\tilde{R}+\tilde{R}^2$, where $\lambda$s are the distribution parameters for each regime and $A$s are the parameters in the transition matrix. $\vartheta^c := (\lambda_1^c,\cdots, \lambda_{\tilde{R}}^c,A_{11}^c, \cdots,A_{\tilde{R}\tilde{R}}^c)$ denotes its unknown true value. 

MSM is able to provide the predictive distribution for future regimes. Given a fixed $\vartheta$, the predictive probability $P(S_{t+1}=l|\xi^t, \vartheta)$ can be calculated using the law of total probability:
\begin{equation}
	\label{predict_dist}
	w_l^t(\vartheta) := P(S_{t+1}=l|\xi^t, \vartheta) = \sum_{k=1}^{\tilde{R}} A_{k,l}P(S_{t}=k|\xi^t, \vartheta). 
\end{equation}
Computing $P(S_{t}|\xi^t, \vartheta)$, known as the filtering problem in MSM, is well-established and can be performed efficiently by applying Bayes rule and law of total probability recursively \citep{fruhwirth2006finite}. 

\subsection{Objective Formulation}\label{sec:obj}
We consider an online simulation optimization problem where input data arrive sequentially and exhibit non-stationarity due to regime switching.
During multiple sequential decision stages, our goal is to identify the optimal decision for each upcoming stage $t+1$. The objective function for stage $t+1$ is formulated as: 
\begin{equation}
	\label{realized-objective}
	\min_{x \in \mathcal{X}} \;  \mathbb{E}_{P(S_{t+1}|\xi^t, \vartheta^c)} \mathbb{E}_{\xi \sim P(\xi_{t+1}|S_{t+1}, \vartheta^c)}[y({x},\xi)].
\end{equation}
Here, $\xi$ is drawn from the distribution corresponding to the regime $S_{t+1}$, with true input parameter $\vartheta^c$. Compared with ordinary stochastic optimization problem \eqref{SO_problem}, the outer layer expectation is derived from $P(S_{t+1}|\xi^t, \vartheta^c)$, the predictive probability for $S_{t+1}$, to account for predictive uncertainty due to potential regime shifts. 
Function \eqref{realized-objective} is built under the best achievable regime prediction, as all parameters governing the regime-switching process (i.e., distribution parameters and regime transition probabilities) are assumed to be known. Hence, it is referred to as the true objective function. However, when making decision, $\vartheta^c$ is unknown and $P(S_{t+1}|\xi^t, \vartheta^c)$ needs to be estimated.  

In practice, we estimate $P(S_{t+1}|\xi^t, \vartheta^c)$ and $\vartheta^c$ with streaming input data. For ease of exposition, we assume that at the beginning of each time stage $t+1$, a single new input data point $\xi_t$ observed from stage $t$ is received. Nonetheless, our approach can be easily extend to situations with a varying number of input data points. To incorporate both uncertainties in regime prediction and parameter estimation, we propose the following Bayesian objective function as an approximation to \eqref{realized-objective}:
\begin{equation}
	\label{BRO-objective}
	\min_{x \in \mathcal{X}} \; \mathbb{E}_{P(\vartheta|\xi^t)} \{\mathbb{E}_{P(S_{t+1}|\xi^t, \vartheta)}[\mathbb{E}_{P(\xi_{t+1}|S_{t+1}, \vartheta)}(y(x, \xi))]\},
\end{equation}
where $\xi^t:=\{ \xi_1, \xi_2,...,\xi_t\}$ denotes the collection of observed input sequence and $P(\vartheta|\xi^t)$ is the posterior of input parameter vector $\vartheta$. Beyond the prediction uncertainty, this Bayesian framework incorporates the input uncertainty through an additional outer layer of expectation, computed with respect to $P(\vartheta|\xi^t)$. Thus, \eqref{BRO-objective} provides a principled way to integrate information from the input data and regime-switching dynamics into the decision process. With MSM, the regime-switching dynamics of input data streams $\xi_t$ can be characterized by a hidden regime Markov chain $S_t$. Given the typically small number of regimes, the three-layer structure can be simplified into a more manageable form using \eqref{predict_dist}:
\begin{equation}
	\label{BRO-objective-simple}
	\mathbb{E}_{P(\vartheta|\xi^t)} \{\mathbb{E}_{P(S_{t+1}|\xi^t, \vartheta)}[\mathbb{E}_{P(\xi_{t+1}|S_{t+1}, \vartheta)}(y(x, \xi))]\} = \mathbb{E}_{P(\vartheta|\xi^t)} \left[\sum_{l=1}^{\tilde{R}} w_l^t(\vartheta) z(x, \lambda_l)\right],
\end{equation}
where $z(x, \lambda_l) := \mathbb{E}_{\xi \sim \tilde{P}(\xi_{t+1}|\lambda_l)}(y(x, \xi))$ is the ordinary stochastic optimization objective with input distribution parameterized by $\lambda_l$.

\section{Theoretical Validation of the Bayesian Approximation}\label{sec:Theory}
During the online process, as the input dataset grows iteratively, the posterior distribution of input parameters is updated, leading to simultaneous improvements in the approximated objective function. This aligns with our intuition that errors from input data estimation should diminish asymptotically. In this section, we theoretically validate this intuition and demonstrate the effectiveness of the Bayesian approximation \eqref{BRO-objective} by establishing the consistency (in Section \ref{Sec: cons_obj}) and asymptotic normality (in Section \ref{Sec: asym_normal}) of its objective function value and optimal solutions. We begin this section by introducing some important notations and definitions.

The hidden regime process $\{S_t, t\geq1\}$ is assumed to be stationary, with initial density $r$ and transition matrix $A$ specified by density $k(s,s^{\prime})$, both with respect to a dominating $\sigma$-finite measure $\mu$. The observation process $\{\xi_t, t\geq1\}$ takes values in $\Xi$. Given $S_i=s$, the conditional distribution of $\xi_i$ has density $g(\cdot|s)$ with respect to a dominating $\sigma$-finite measure $\nu$ on $\Xi$. The joint Markov transition kernel $Q$ is dominated by $\mu\otimes\nu$, with density $q(z,\cdot)$. $\vartheta$ takes values in the space $\Theta$, equipped with Borel $\sigma$-algebra $ \mathcal{B}_{\Theta}$. Let $\mathcal{B}_{\Xi}^N$ denote the Borel $\sigma$-algebra on the space of all infinite sequences in $\Xi$. For fixed $\vartheta$, let $P_\vartheta^t$ and $P_\vartheta^N$ denote the laws of $\xi^t$ and infinite observation sequence, respectively. The probability law under $\vartheta^c$ is denoted by $P^{\vartheta^c}$, with expectation $\mathbb{E}^{\vartheta^c}$. The distance metrics $\mathbb{D}$ and $\Vert \cdot \Vert_\textup{TV}$ are defined according to Definitions 3.5 and 4.1 in \citet{wu2018bayesian}. The former measures the distance between two sets (e.g., of optimal solutions), while the latter quantifies the difference between probability measures (e.g., in the Central Limit Theorem (CLT)). For notational simplicity, we denote the posterior $P(\vartheta|\xi^t)$ by $P_t$, and define 
$\widetilde{K}(x,\vartheta):=\mathbb{E}_{P(S_{t+1}|\xi^{t},\vartheta)}[\mathbb{E}_{P(\xi_{t+1}|S_{t+1},\vartheta)}(y(x,\xi_{t+1}))]$,
where $\widetilde{K}(x,\vartheta^c)$ is indeed the true objective \eqref{realized-objective}.

\begin{defn}[Density function of observations]
	The measure $P_\vartheta^t$ admits the following density of $\xi^t$ with respect to the product measure $\nu^t$:
	$$p_{\vartheta}(\xi^t):= \int_{S^t} r_{\vartheta}(s_1)\prod_{j=1}^{t-1}k_{\vartheta}(s_j,s_{j+1})\prod_{i=1}^t g_{\vartheta}(\xi_i|s_i)\prod_{i=1}^t\mu(ds_i).$$
\end{defn}

\begin{defn}[Posterior distribution of $\vartheta$]
	Given $\xi^t$ and a prior distribution $P_0$, for any measurable set $C\in\mathcal{B}_{\Theta}$, the posterior distribution of $\vartheta$ on $(\Theta, \mathcal{B}_{\Theta})$ is given by:
	$$P_t(C):= \frac{\int_C p_{\vartheta}(\xi^t)P_0(d\vartheta)}{\int_{\Theta} p_{\vartheta}(\xi^t)P_0(d\vartheta)}.$$ 
\end{defn}

\begin{defn}[Joint measure]
	Let $\eta$ be a measure defined by: $\eta(A_1\times A_2) := \int_{A_1} \int_{A_2} P_\vartheta^N(d\xi)P_0(d\vartheta)$ for $\forall A_1\in \mathcal{B}_{\Theta}, A_2\in \mathcal{B}_{\Xi}^N.$ 
\end{defn}

\subsection{Consistency of Objectives}
\label{Sec: cons_obj}
In this subsection, we first derive the posterior consistency for MSM parameters in Lemma \ref{cons_msm}. Next, we prove the pointwise consistency of our objective function  
and the consistency of optimal solutions and values in Theorems \ref{consistency_obj} and \ref{theorem2}, respectively. The proofs are provided in Appendix C.

\begin{assumption}
	\label{assum1}
	\begin{enumerate}[(i)]
		\item $\mathbb{E}_{\tilde{P}(\xi_{t+1}|\vartheta)}[y(x,\xi_{t+1})]$ is continuous on $\Theta$ for every fixed x.
		\item The emission distribution $\tilde{P} (\xi_t|\lambda_i)$ is continuous in $\vartheta\in\Theta$ for $i=1,\cdots,\tilde{R}$. 
	\end{enumerate}
\end{assumption}

Assumption \ref{assum1} is readily satisfied when $\tilde{P}(\cdot \mid \vartheta)$ belongs to a regular exponential family, such as exponential distribution and Gaussian distribution. In such case, the density function varies smoothly with respect to the natural parameters, and $\mathbb{E}[y(x,\xi)]$ is continuous in $\vartheta$ as long as $y(x,\xi)$ is integrable with respect to $\tilde{P}(\cdot \mid \vartheta)$.

\begin{assumption}
	\label{assum2}
	\begin{enumerate}
		[(i)]
		\item For all $\vartheta \in \Theta$, and $z_1=(s_1,\xi_1)$, $z_2=(s_2, \xi_2)\in Z^{\prime}=S\times \Xi$, the transition density satisfies $q_{\vartheta}(z_1,z_2)=k_{\vartheta}(s_1,s_2)g_{\vartheta}(s_2,\xi_2)$, with $q_{\vartheta}(z_1,z_2)>0\quad$ a.s. $(\mu\otimes\nu)$. 
		\item $P_0(\Theta)<\infty$, and $P_0$ is proper, i.e. for all $\delta>0$, $P_0(\{\vartheta \in \Theta: \Delta(\vartheta^c,\vartheta)\leq\delta\})>0$, where 
		$\Delta(\vartheta^c,\vartheta):= \liminf_{n\rightarrow\infty} \frac{1}{n} \sum_{l=1}^{n-1}\log \frac{q_{\vartheta^c} (z_l, z_{l+1})}{q_{\vartheta} (z_l, z_{l+1})}$.
	\end{enumerate}
\end{assumption}

Assumption 2 (i) implies that the transition kernel density admits a factorized form and is strictly positive almost surely. Assumption 2 (ii) ensures the prior distribution $P_0$ does not concentrate around parameters whose likelihood is too small asymptotically. Commonly used priors, such as uniform distribution and gamma distribution, typically satisfy this condition. Under these mild assumptions, we first prove the following Lemma \ref{cons_msm} to obtain the strong consistency of posterior distribution for input parameters in MSM with a finite regime space.

\begin{lem}[Consistency of posterior distribution for MSM]
	\label{cons_msm}
	Suppose Assumptions (B1)-(B6) in \citet{de2008asymptotic}, and Assumption \ref{assum2} hold true, then the posterior distribution $P_t\Rightarrow \delta_{\vartheta^c}$ a.s.$(P_{\vartheta^c}^N)$, where $\delta_{\vartheta^c}$ is a point mass on true parameter $\vartheta^c$, i.e. 
	$P_{\vartheta^c}^N(P_t \stackrel{t\rightarrow\infty}{\Rightarrow}\delta_{\vartheta^c})=1$.    
\end{lem}

This result is specifically built upon the general consistency analysis of MSMs with measurable regime spaces by \citet{douc2020posterior}. This lemma plays a crucial rule in establishing the consistency of our proposed objective. Based on Lemma \ref{cons_msm}, we establish the following Theorem \ref{consistency_obj} for consistency of the function value. 

\begin{thm}[Consistency of the objective function]
	\label{consistency_obj}
	Under Assumptions (B1)-(B6) in \citet{de2008asymptotic}, Assumption \ref{assum1} and \ref{assum2}, as $t \rightarrow \infty$, for every fixed $x\in \mathcal{X}$, the following holds a.s. $(P_{\vartheta^c}^N):$ 
	$$\mathbb{E}_{P_t}\{{\mathbb{E}_{P(S_{t+1}|\xi^{t},\vartheta)}[\mathbb{E}_{P(\xi_{t+1}|S_{t+1},\vartheta)}(y(x,\xi_{t+1}))]}\}\rightarrow{\mathbb{E}_{P(S_{t+1}|\xi^{t},\vartheta^c)}[\mathbb{E}_{P(\xi_{t+1}|S_{t+1},\vartheta^c)}(y(x,\xi_{t+1}))]}.$$
\end{thm}

Theorem \ref{consistency_obj} shows that as data accumulates, our proposed objective value converges to the true objective value $\widetilde{K}(x, \vartheta^c)$.  

\begin{thm}[Consistency of the optimal solutions and optimal values]
	\label{theorem2}
	Suppose Assumption 3.2 in \citet{wu2018bayesian} holds, $\Theta$ and $\mathcal{X}$ are compact, and $\widetilde{K}(x,\cdot)$ is continuous on $\Theta$ for every $x\in\mathcal{X}$. Let $M_t := \arg \min_{x\in\mathcal{X}} \mathbb{E}_{P_t}[\widetilde{K}(x,\vartheta)]$ be the set of optimal solutions, and $M := \arg\min_{x\in\mathcal{X}} \widetilde{K}(x,\vartheta^c)$ be the set of true optimal solutions. If there exists a measurable function $\psi: \Theta\rightarrow \mathbb{R}^{+}$ with $\left|\widetilde{K}(x,\vartheta)-\widetilde{K}(y,\vartheta)\right|\leq \psi(\vartheta) \Vert x-y \Vert$, $\forall x,y\in\mathcal{X}$ and $\int_{\Theta} \psi(\vartheta)\lambda(d\vartheta)<\infty$, then as $t\rightarrow\infty$,
	$$\mathbb{D}(M_t, M)\rightarrow 0 \quad and \quad \min_{x\in\mathcal{X}} \; \mathbb{E}_{P_t}[\widetilde{K}(x,\vartheta)]\rightarrow \min_{x\in\mathcal{X}}\; \widetilde{K}(x,\vartheta^c) \quad a.s. (\eta).$$
\end{thm}

These theorems validates that our proposed objective function provides a reliable and consistent approximation asymptotically.

\subsection{Asymptotic Normality of Objectives}
\label{Sec: asym_normal}
We further explore the asymptotic normality of the function value, optimal values, and optimal solutions of the proposed approximation. Asymptotic normality results characterize the distributional behavior of the solutions as $t$ grows. These results enable us to assess the rate of convergence and support statistical inferences. For example, they allow for the construction of confidence intervals for the optimal values or solutions, quantifying the uncertainty inherent in using the approximation. Moreover, by specifying a tolerance level for the variance (or the width of the confidence interval), one can determine the minimum amount of real-world data required to achieve the target variance using estimators of the asymptotic variances. In this section, we demonstrate that, asymptotically, $\sqrt{t}$ times the difference between our estimated optimal values/solutions and the true values under perfect distribution knowledge follows a normal distribution. We first establish a Bayesian CLT for MSM in terms of total variation in Lemma \ref{lemma-added}.
\begin{lem}[Bayesian total variation CLT for MSM]
	\label{lemma-added}
	Under Assumptions (B1)-(B6) and (A5) in \citet{de2008asymptotic}, as $t\rightarrow \infty$, we have $\Vert P_{Z_t}-\mathcal{N}(\Delta_t, J_c^{-1})\Vert_\textup{TV}\rightarrow 0$ in probability $(P_{\vartheta^c}^N)$ as $t\rightarrow \infty$, where $Z_t (\vartheta) := \sqrt{t}(\vartheta-\vartheta^c)$ and $P_{Z_t} := P_t \circ Z_t^{-1}$. $J_c$ is the non-singular Fisher information matrix evaluated at true parameter $\vartheta^c$, and $\Delta_t:=\frac{1}{\sqrt{t}}J_c^{-1}\nabla L_t(\vartheta^c) \Rightarrow \mathcal{N}(0,J_c^{-1})$, where $L$ is the log-likelihood function. 
\end{lem}

While \citet{de2008asymptotic} establish total variation convergence for the posterior centered at a data-dependent point involving $\Delta_t$, our result demonstrates total variation convergence when centering directly at $\vartheta^c$. This refinement enables a more direct and interpretable asymptotic characterization, and is central to establishing the asymptotic normality of the objective functions.

\begin{thm}[Asymptotic normality at a fixed ${x}$]
	\label{mean-normality}
	Under Assumption 4.1 in \citep{wu2018bayesian}, Assumptions (B1)-(B6) and (A5)-(A6) in \citep{de2008asymptotic}, if $\widetilde{K}$ is continuous on $\Theta$ and differentiable at $\vartheta^c$, as $t \rightarrow \infty$, we have $\sqrt{t}\{\mathbb{E}_{P_t}[\widetilde{K}(\vartheta)]-\widetilde{K}(\vartheta^c)\}\Rightarrow \mathcal{N}(0,\sigma_x^2),$ where $\sigma_x^2:=\nabla_\vartheta \widetilde{K}(x,\vartheta^c)^T J_c^{-1} \nabla_\vartheta \widetilde{K}(x,\vartheta^c)^T.$ 
\end{thm}

\begin{thm}[Asymptotic normality of optimal values]
	\label{opt-normality}
	Suppose that Assumption \ref{assum2}, Assumptions (B1)$\sim$(B6) in \citep{de2008asymptotic}, and Assumption 4.1 in \citet{wu2018bayesian} all hold. Also assume for $\vartheta \sim P_t$ that $\sqrt{t}(\vartheta-\vartheta^c)$ and $\sqrt{t}[\widetilde{K}(\vartheta)-\widetilde{K}(\vartheta^c)]$ have positive densities for all $t$ a.s. $(P_{\vartheta^c}^N)$. Further suppose $\mathcal{X}$ is a compact set, $\widetilde{K}$ is continuous on $\mathcal{X}\times \Theta$, and $\widetilde{K}(x,\cdot)$ is differentiable at $\vartheta^c$ for all $x\in \mathcal{X}$, where $\nabla_{\vartheta}\widetilde{K}(\cdot, \vartheta^c)$ is continuous on $\mathcal{X}$. Then,
	$$\sqrt{t}(\min_{x\in\mathcal{X}}\;\mathbb{E}_{P_t}[\widetilde{K}(x,\vartheta)]-\min_{x\in\mathcal{X}}\;\widetilde{K}(x,\vartheta^c))\Rightarrow \min_{x\in S^*}\; Y_x,$$
	where $S^* := \arg\min_{x\in \mathcal{X}}\widetilde{K}(x,\vartheta^c)$ and $Y_x := \nabla_\vartheta \widetilde{K}(x,\vartheta^c)^T Z^*$, where $Z^* \sim \mathcal{N}(0, J_c^{-1}).$
\end{thm}

Both consistency and asymptotic normality provide theoretical support for the claim that the proposed approximated objective function serves a well-founded approximation to the original objective function under conditions of regime switching.

\section{Regime-Switching Online Bayesian Simulation Optimization Algorithm}\label{sec:Algorithm}

In this section, we develop a metamodel-based algorithm to solve the proposed optimization problem \eqref{BRO-objective}. Specifically, in Section \ref{Sec: GP}, we present a Gaussian process metamodel for the objective function, which integrates joint modeling of the decision vector $x$ and the distribution parameter $\lambda$ over the posterior and regime predictive distributions. This joint model is critical for leveraging simulation results from all prior stages. In Section \ref{Sec: GP Algo}, we develop a regime-aware expected improvement function to sequentially select design points and input parameters for experimental runs. In Section \ref{sec:Algorithm_convergence}, we provide the convergence analysis.

\subsection{The GP Model for the Objective Function}
\label{Sec: GP}
Denote $\tilde{g}_t (x)$ as the Bayesian approximation \eqref{BRO-objective} that we are solving in stage $t$: $\tilde{g}_t (x) = \mathbb{E}_{P(\vartheta|\xi^t)} \left[\sum_{l=1}^{\tilde{R}} w_l^t(\vartheta) z(x, \lambda_l)\right]$, where all parameters in $\vartheta$ are unknown. For the proposed algorithm, we need to build a metamodel for $\tilde{g}_t (x)$. Meanwhile, the simulation results from past stages should be used to build the model. It is important to note that these results are run under different input distributions at stage $t$, driven by variations in the input dataset and the regime-shift dynamics. Consequently, they cannot be directly used to estimate the response at any point $x$, even if that point was evaluated in previous stages. However, discarding these results entirely is inefficient, as they contain valuable information, particularly from more recent stages when the dataset and input modeling closely resemble the current stage. Moreover, some simulation experiments can be time-intensive, and the number of experiments feasible in each stage may be limited. In such cases, leveraging informative results from past stages, rather than restarting the optimization process from scratch, is critical. This represents a key research focus in online simulation optimization \citep{wu2024data}. To construct a model for $\tilde{g}_t (x)$ that leverages past results, we first develop a model $Z(x, \lambda)$ for $z(x, \lambda)$. This model explicitly accounts for variations in the input distributions, through including $\lambda$ as the model input, of previous experiment results. By fitting it with all past results, the model enables information sharing to infer the response value at the current stage for any desired $x$ and $\lambda$. We then derive a model $\hat{G}_t (x)$ for $\tilde{g}_t (x)$ through integrations of $Z(x, \lambda)$ over the regime predictive distribution and the input posterior.

We first fit a stochastic GP model $Z(x, \lambda)$ for $z(x, \lambda)$ with respect to both $x$ and $\lambda$. $Z(x, \lambda)$ is assumed to be a zero mean GP with the covariance structure $\text{Cov}[ Z(x_j, \lambda_j), Z(x', \lambda')] = \sigma_g^2 K((x_j, \lambda_j),(x', \lambda'))$, where $\sigma_g$ is the spatial variance of the GP model and $K((x_j, \lambda_j),(x', \lambda'))$ is a kernel function. Popular choices of $K$ include the Gaussian kernel and the Mat\'{e}rn kernel. In this work, we use the Gaussian kernel:
$K((x_j, \lambda_j),(x', \lambda')) = \exp(-d^2[(x_j, \lambda_j), (x', \lambda')] ),$
where $d^2[(x, \lambda), (x', \lambda')]=  \sum_{i=1}^{d_x}
\frac{(x_i-x_i^{\prime})^{2}}{2\theta_{1,i}^2} + \sum_{j=1}^{d_\lambda}\frac{(\lambda_i-\lambda_i^{\prime})^{2}}{2\theta_{2,j}^2}$. $\theta$s are the length scale parameters of the Gaussian kernel, and $d_x$ and $d_\lambda$ are the dimensions of $x$ and $\lambda$, respectively. In this model, the random simulation result obtained under $x$ and input parameter $\lambda$ is assumed to be $y(x, \lambda) = z(x, \lambda) + \epsilon$, where $\epsilon$ is a normal noise following $\mathcal{N}(0, \sigma_\epsilon^2)$ independent of the GP.   

Suppose the combinations of $x$ and $\lambda$ to run the experiments in all past stages are $\{(x_1, \lambda_1), \cdots, (x_n, \lambda_n)\}$, each with $m$ replications. Let $Y_n=[(\bar{y}(x_1, \lambda_1), \cdots,\bar{y}(x_n,\lambda_n)]^T$ denote the vector of averaged simulation output. Conditional on $Y_n$, the posterior distribution of $Z(x, \lambda)$ is given by $Z_n(x, \lambda) := Z|\bar{Y_n}\sim \text{GP}(m_n, k_n)$, where the posterior mean and covariance are 
$$m_n(x, \lambda) = c(x,\lambda)^T R^{-1}Y_n,$$
$$k_n((x,\lambda), (x^{\prime},\lambda^{\prime})) = c((x,\lambda), (x^{\prime},\lambda^{\prime})) -{ c({x, \bf \lambda})}^T{R}^{-1}{{ c}({x^{\prime}, \lambda^{\prime} })}.$$ 
Here, $c(x,\lambda)=[\text{Cov}[Z(x,\lambda), Z(x_1, \lambda_1)], \cdots, \text{Cov}[Z(x, \lambda), Z(x_n, \lambda_n)]]^T$ is an $n\times1$ covariance vector between $(x, \lambda)$ and past points. $R = R_z + R_\epsilon$ is the $n\times n$ covariance matrix, where $R_z$ is the covariance matrix among existing design points from the GP and $R_\epsilon$ is the diagonal noise covariance matrix of the simulation results.

Based on model $Z$, we can build the model for $\tilde{g}_t (x)$ through numerical integration for the two outer layers of expectation in \eqref{BRO-objective}. For the outermost layer, as $\vartheta$ includes both distribution parameters and transition probabilities, the posterior distribution $P(\vartheta|\xi^t)$ is analytically intractable. Therefore, we employ Markov Chain Monte Carlo (MCMC) techniques to generate posterior samples and approximate the outer expectation with sample average approximation. The sampling procedures consist of the following steps. First, we employ the forward-filtering backward-sampling method from  \citet{fruhwirth2006finite} to sample the hidden regime sequence: 
$$S_T \sim P(S_T \mid \xi_{1:T}, \vartheta), \; S_t \mid S_{t+1} \sim P( S_{t+1} \mid S_t, \vartheta ) P( S_{t+1} \mid \xi_{1:T}, \vartheta).$$
Next, we impose a symmetric Dirichlet prior $\text{Dir}(1, \cdots, 1)$ on each row of the transition matrix $A_i = (A_{i,j})_{j=1}^{\tilde{R}}$, serving as a conjugate prior for the multinomial distribution. With the count $n_{i,j}$ representing the number of transitions from regime $i$ to regime $j$ observed in the sampled regime sequence, the posterior of $A_i$ becomes $\text{Dir}(1+n_{i,1}, \cdots, 1+n_{i,\tilde{R}})$. Finally, for the distribution parameters, we use the No-U-Turn Sampler to generate posterior samples \citep{hoffman2014no}, with weakly-informative priors. Each of the generated samples can then be used in MSM to compute the weights $w$ in \eqref{predict_dist} to calculate the middle layer of expectation. Therefore, we can derive the following GP model $\hat{G}_t(x)$ for $\tilde{g}_t (x)$ (here, $n$ denotes the total number of design points accumulated from previous $t-1$ time stages and initial designs):
\begin{equation}
	\label{weighted_gp}
	\hat{G}_t(x) = \frac{1}{N_{MC}}\sum_{i=1}^{N_{MC}} \sum_{l=1}^{\tilde{R}}w_l^{(t, i)} Z_n(x, {\lambda}_l^{(t, i)}),
\end{equation} 
where the mean and covariance are given by:
$$\mathbb{E}[\hat{G}_t(x)] = \mu_{t}(x) = \frac{1}{N_{MC}}\sum_{i=1}^{N_{MC}} \sum_{l=1}^{\tilde{R}} w_l^{(t, i)} m_n(x, \lambda_l^{(t,i)}),$$
\begin{equation*}
	\text{Cov}[\hat{G}_t(x), \hat{G}_t(x')] = v_t(x, x') = \frac{1}{N_{MC}^2} \sum_{i=1}^{N_{MC}} \sum_{j=1}^{N_{MC}} \sum_{l=1}^{\tilde{R}} \sum_{m=1}^{\tilde{R}} 
	w_l^{(t,i)} w_m^{(t,j)} k_n\left((x, \lambda_l^{(t,i)}), (x', \lambda_m^{(t,j)})\right). 
\end{equation*}

The GP model, despite its effectiveness, can be computationally restrictive in online settings due to its $O(n^3)$ complexity. To enhance computational efficiency when processing large volumes of simulation results, we can adopt approximate GP models such as streaming sparse GP approximations \citep{bui2017streaming}, or a moving window approach to discard outdated simulation outputs \citep{wu2024data} to maintain manageable computational demands.

\subsection{The GP-Based Optimization Algorithm}
\label{Sec: GP Algo}
At a fixed time stage $t$, we aim to identify the minimizer of the objective function \eqref{BRO-objective} with the help of the GP model. To guide the selection of the next design point, we develop a regime-aware Expected Improvement (EI) function as the searching criterion. EI function is a widely-applied acquisition function for GP-based algorithm (also termed as Bayesian optimization) since \citet{jones1998efficient}, due to its ability to balance exploration and exploitation. Following the approach in \citet{pearce2017bayesian}, we assume a new hypothetical sample $(x_{n+e_t+1}, \lambda_{n+e_t+1})$ is generated ($e_t$ represent the current number of design points collected during the $t$-th time stage), and then get the predictive distribution of $\hat{\mu}_{n+e_t+1}$ based on the hypothetical evaluations at $(x_{n+e_t+1}, \lambda_{n+e_t+1})$:
\begin{equation*}
	\label{predicted_mu_weight}
	\hat{\mu}_{n+e_t+1}(x|x_{n+e_t+1}, \lambda_{n+e_t+1}) \sim \mathcal{N} \left(\mu_{n+e_t}(x), \widetilde{\sigma}^2(x|x_{n+e_t+1}, \lambda_{n+e_t+1})\right),
\end{equation*}
where the variance is given by:
$$\widetilde{\sigma}^2(x|x_{n+e_t+1}, \lambda_{n+e_t+1}) = \left(\frac{1}{N_{MC}}\sum_{i=1}^{N_{MC}} \sum_{l=1}^{\tilde{R}}w_l^{(t, i)} \frac{k_{n+e_t}((x, \lambda_l^{(t,i)}), (x_{n+e_t+1}, \lambda_{n+e_t+1}))}{\sqrt{k_{n+e_t}((x_{n+e_t+1}, \lambda_{n+e_t+1}), (x_{n+e_t+1}, \lambda_{n+e_t+1}))+\frac{\sigma_{\epsilon}^2}{m}}}\right)^2.$$
The EI acquisition function in our setting is defined as:
\begin{equation*}
	\label{EI_joint}
	\text{EI}(x, \lambda) = \mathbb{E}[(T-\hat{\mu}_{n+e_t+1}(x|x,\lambda))^+] = \Delta\Phi(\frac{\Delta}{\widetilde{\sigma}(x)}) + \widetilde{\sigma}(x)\phi(\frac{\Delta}{\widetilde{\sigma}(x)}),
\end{equation*}
where the current best value $T=\min\{\mu_{n+e_t}(x_1),\cdots,\mu_{n+e_t}(x_{{n+e_t}})\}$, which is the minimum posterior mean after the previous overall $n+e_t$ design points. $\Delta=T-\mu_{n+e_t}(x)$, and $\widetilde{\sigma}(x) = \widetilde{\sigma}(x|x,\lambda)$. $\Phi$ and $\phi$ are the cdf and pdf of the standard normal random variable. The next design point $(x_{n+e_t+1}, \lambda_{n+e_t+1})$ is chosen with $\arg\max_{x, \lambda}\text{EI}(x, \lambda)$ to achieve the greatest expected improvement.

\begin{algorithm}
	\caption{The RSOBSO Algorithm}
	\label{algo-BRO}
	\begin{algorithmic}[1]
		\State \textbf{Initialization:} Generate $\{x_1,\cdots,x_{n_0}\}$ using LHS; Given initial input data $\xi^0$, and prior $P_0(\vartheta)$, generate $n_0$ sets of samples from posterior $P(\vartheta|\xi^0)$ to obtain initial points $\{(x_1, {\lambda_1^{(1)}}),\cdots,(x_{n_0}, {\lambda_1^{(n_0)}}), \cdots, (x_1, {\lambda_{\tilde{R}}^{(1)}}),\cdots,(x_{n_0}, {\lambda_{\tilde{R}}^{(n_0)}})\}$; Set $n\leftarrow n_0 * \tilde{R}$;
		\State Run simulations at initial points with $m$ replications each and get observation vector ${\bar{Y}^{(0)}_n}$; 
		\For{$t = 0$ to $t_{max}-1$}
		\For{$e_t = 0$ to specified number of optimization steps $u-1$}
		\State Construct GP models $Z^{(t)}_n({x}, \lambda)$ and $\hat{G}_t({x})$ based on $\bar{Y}^{(t)}_n$ and $P(\vartheta|\xi^{t})$;
		\State Select next evaluation point $({x}^{(t)}_{n+e_t+1}, \lambda^{(t, n+e_t+1)})$ using the regime-aware EI function;
		\State Run simulation experiments at $({x}^{(t)}_{n+e_t+1}, \lambda^{(t, n+e_t+1)})$ with $m$ replications and obtain \text{\hspace{2.8em}} $\bar{y}({x}^{(t)}_{n+e_t+1}, \lambda^{(t, n+e_t+1)})$, set ${\bar{Y}^{(t)}_{n+e_t+1}} \leftarrow [ {\bar{Y}_{n}^{(t)}},\bar{y}({x}^{(t)}_{n+e_t+1}, \lambda^{(t, n+e_t+1)})]^T$ and $n\leftarrow n+e_t+1$; 
		\EndFor
		\State \textbf{return} $\hat{x}_t \leftarrow \arg\min_{x}\mu_t({x})$;
		\State Receive new input data $\xi_{t+1}$ and update $P(\vartheta|\xi^{t+1})$ and set $\bar{Y}_{n}^{(t+1)} \leftarrow \bar{Y}_{n}^{(t)}$;
		\EndFor
	\end{algorithmic}
\end{algorithm}

The RSOBSO algorithm is outlined in Algorithm \ref{algo-BRO}. Steps 1 and 2 generate initial samples to run simulations. From step 3 onward, the algorithm enters a loop over time stage $t$. Steps 4 to 8 iteratively select design point using the regime-aware EI function and evaluate it. Step 9 yields the optimal decision $\hat{x}_t$ for stage $t$. Steps 10 updates the posterior of $\vartheta$ upon receiving new data $\xi_{t+1}$ and carries forward all simulation results.  

\subsection{Algorithm Convergence Analysis}\label{sec:Algorithm_convergence}
In this section, we establish the convergence analysis of our algorithm. For ease of exposition, we assume that the solutions to the objective functions are all unique. Denote $g_t(x)$ as the true objective \eqref{realized-objective}. Recall that $\tilde{g}_t(x)$ is the proposed Bayesian approximation and $\hat{x}_t$ is the solution returned by Algorithm 1. We further denote $x_t^* = \arg \min g_t(x)$ and $\tilde{x}_t^* = \arg \min \tilde{g}_t(x)$. The analysis can be divided in two steps. In the first step, we show that for any time stage $t$, the value of the optimal solution returned by Algorithm 1, $\tilde{g}_t(\hat{x}_t)$, will converge to the optimal value of the objective function in this stage $\tilde{g}_t(\tilde{x}_t^*)$, given unlimited computing resources. In the second step, we show that when both the decision stage and the computing resources tend to infinity, the value of the returned solution, $g_t(\hat{x}_t)$, will converge to the true optimal value $g_t(x_t^*)$. Similar to \cite{wang2022multilevel} and \cite{wang2020nonparametric}, the analysis is based on the following two assumptions:
\begin{assumption}
	\label{ass_desigh}
	The decision space $\mathcal{X}$ is compact.
\end{assumption}
\begin{assumption}
	\label{ass_para}
	The hyperparameters of the GP model are known.
\end{assumption}

The two steps can be summarized into the two following theorems.

\begin{thm}
	\label{convergence_singlestage}
	In decision stage $t$, when $u, m, N_{MC}$ tend to infinity, we have $\tilde{g}_t(\hat{x}_t) \rightarrow_p \tilde{g}_t(\tilde{x}_t^*)$, under Assumptions \ref{ass_desigh} and \ref{ass_para}.
\end{thm}
To show this theorem, we first prove that the evaluated points selected by EI is everywhere dense, according to the analysis of \cite{locatelli1997bayesian}. Then, we can show that the prediction $\mu_t(x)$ converges to $\tilde{g}_t(x)$. Following the similar reasoning to \cite{wang2020nonparametric}, we achieve the convergence results. 

\begin{thm}
	\label{convergence_whole}
	When $t, u, m, N_{MC}$ tend to infinity, we have $g_t(\hat{x}_t) \rightarrow_p g_t(x_t^*)$, under Assumptions \ref{ass_desigh} and \ref{ass_para}.
\end{thm}
This theorem is a result of both Theorem \ref{convergence_singlestage} and the theoretical analysis in Section \ref{sec:Theory}, combining the convergence of the algorithm in a single decision stage and the convergence of the approximated objective function $\tilde{g}_t(x)$ across different stages. The detailed proof are provided in Appendix F.  

\section{Numerical Experiment}\label{sec:Experiment}
We evaluate the empirical performance of RSOBSO algorithm on two synthetic benchmark problems in Section~\ref{Sec: synthetic} and two real-world applications: inventory management (Section~\ref{Sec: inventory_problem}) and portfolio optimization (Section~\ref{Sec: portfolio}). To make the experiments comprehensive and results representative, we compare RSOBSO with four commonly used methods:
\begin{enumerate}[(i)]
	\item \textbf{RSOPSO} (Regime-Switching Online Plug-in Simulation Optimization): a regime-switching plug-in method that uses the posterior mean of MSM parameters, ignoring input uncertainty.
	\item \textbf{NOBSO} (Non-Regime-Switching Online Bayesian Simulation Optimization): a non-regime-switching Bayesian method, ignoring prediction uncertainty.
	\item \textbf{NOPSO} (Non-Regime-Switching Online Plug-in Simulation Optimization): a non-regime-switching plug-in variant that uses the posterior mean, ignoring both prediction and input uncertainties.
	\item \textbf{NOKSO} (Non-Regime-Switching Online KDE Simulation Optimization): a non-regime-switching nonparametric method that models the input distribution via kernel density estimation (KDE) and constructs a GP directly on $Z(x)$.
\end{enumerate}
The corresponding optimization objective functions of these methods are detailed in Appendix E.

To ensure a fair comparison, all algorithms use the same evaluation budget. Based on historical observations $\xi$ over initial $h$ time stages, we optimize the problems over additional $t_{max}$ time stages. For each time stage $t$, we can estimate the optimal solution $\hat{x}_t$ and denote the true minimizer at period $t$ as $x^*_{S_t}$. The GAP value for time stage $t$ is $z(\hat{x}_t, \lambda_{S_t}) - z(x^*_{S_t}, \lambda_{S_t})$, and the cumulative GAP till period $t$ is $\sum_{i=h+1}^{h+t}\left[z(\hat{x}_i, \lambda_{S_i}) - z(x^*_{S_i}, \lambda_{S_i})\right]$, which is our performance measure adopted. The GAP measures the difference in the objective under the true `realized' regime on the next day between the estimated and true optimal solution. This metric effectively captures both regime prediction and input parameter uncertainties. Detailed experiment settings for this section, such as the prior distributions of parameters and  the decision spaces, are provided in Appendix D. 

\subsection{Synthetic Test Problems}
\label{Sec: synthetic}
We construct two synthetic test problems to illustrate the performance of RSOBSO. The regime data in this synthetic setting are generated from the known transition matrix, as detailed in Appendix D. Nevertheless, during the numerical experiments, the transition probabilities are treated as unknown parameters and estimated within the algorithm.

\subsubsection{4-Regime Exponential Emissions Case}
This is a univariate quadratic problem adapted from \citet{wu2024data}: $y(x, \xi) = (x-\xi)^2 + 10\xi$, where $\xi$ follows the exponential distributions with regime-switching rate $\lambda_i$. The rates for the four regimes are $\lambda_i \in \{1/30, 1/20, 1/10, 1\}$, with optimal solution $x^{*,i}=1/\lambda_i$ and value $y^{*,i}=10/\lambda_i+1/\lambda_i^2$. 

\subsubsection{3-Regime Gaussian Emissions Case}
\begin{figure}[h]
	\centering
	\subfloat{\includegraphics[width=0.46\textwidth]{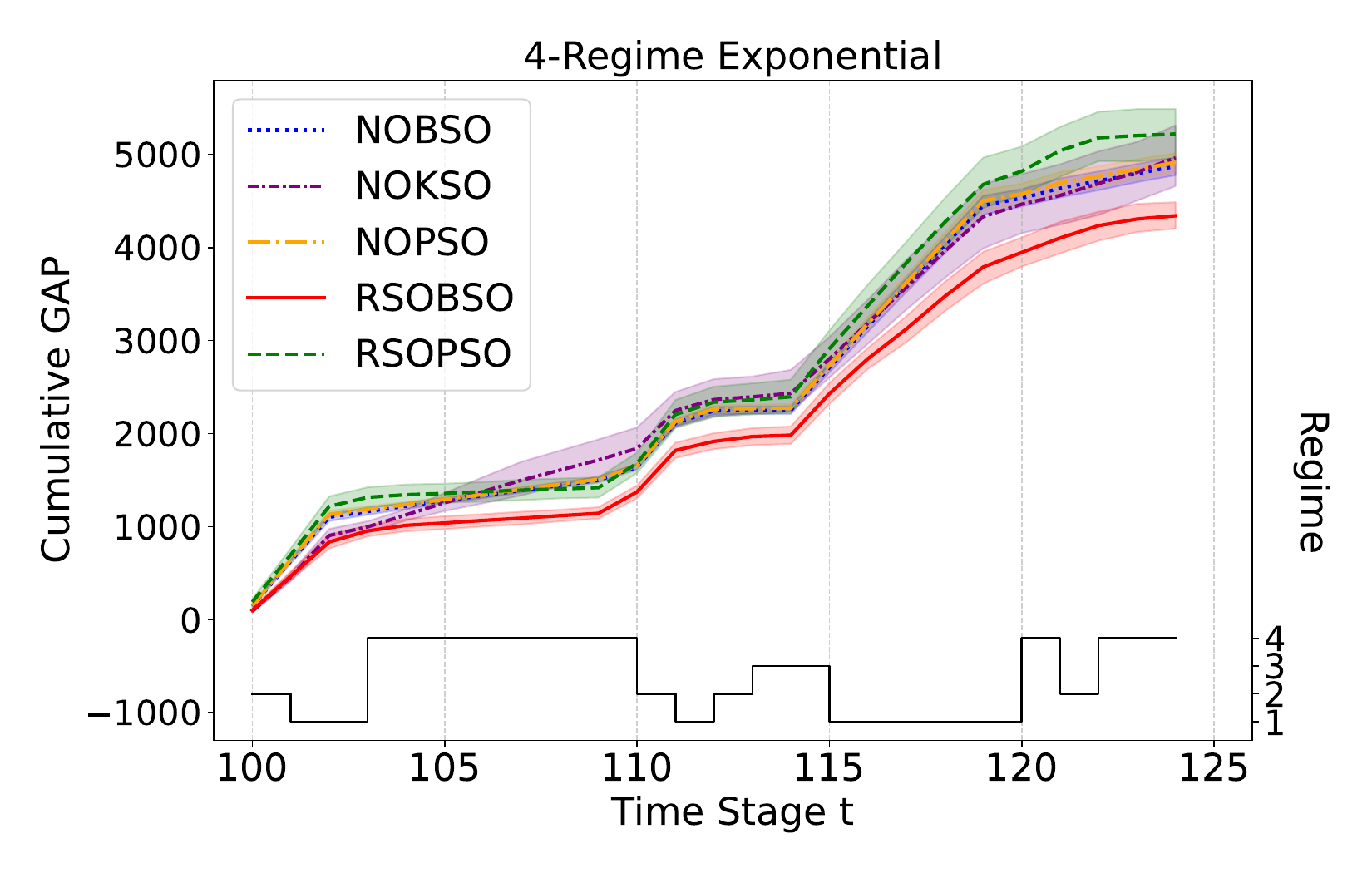}}
	\hfill
	\subfloat{\includegraphics[width=0.46\textwidth]{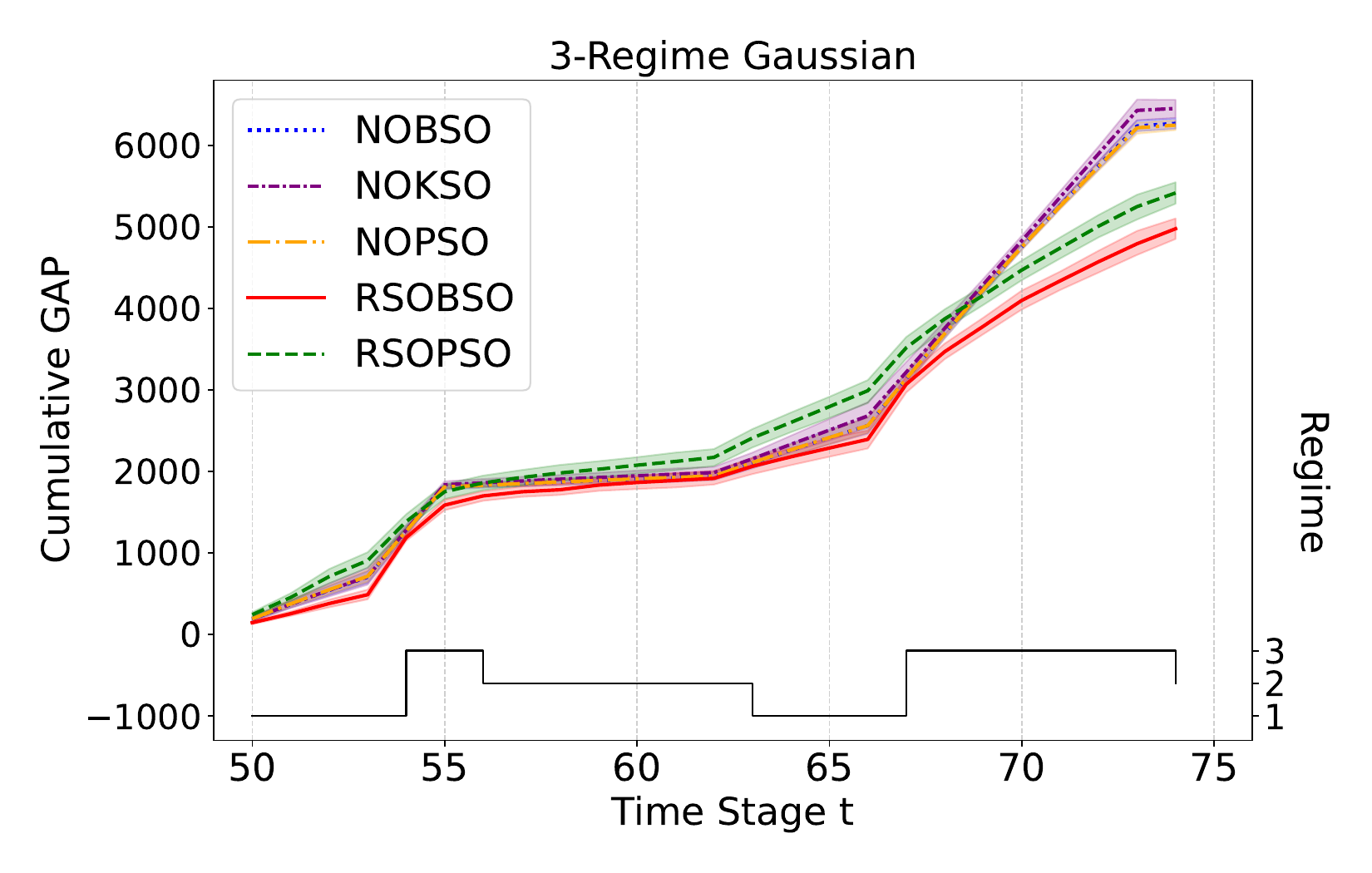}}
	\caption{Cumulative GAP for synthetic test problems (black line indicates regimes).}
	\label{fig:combined_gap}
\end{figure}

This is a multivariate quadratic problem adapted from \citet{liu2024bayesian}: $y(x, \xi) = (x_1 - 10)^2 + (x_2  -20)^2+\xi(4x_1+8x_2)$, where $\xi$ follows $\mathcal{N}(\mu_i, \sigma^2)$, and $\sigma = 3$ (known). The means for the three regimes are $\mu_i \in \{2, 4, 10\}$, with optimal solution $x^{*,i}=(10-2\mu_i, 20-4\mu_i)$, and value $y^{*,i}=200\mu_i-20\mu_i^2$. 
We plot the regimes and the cumulative GAP of the two problems in Figure \ref{fig:combined_gap}. 

\subsubsection{Analysis of Synthetic Test Problems}
We analyze the synthetic test problems with a focus on input distribution estimation and regime-switching handling.

NOBSO and NOPSO assume a single regime and produce similar parameter estimates, often achieving low GAP in one regime but high in others due to their inability to capture regime dynamics. NOKSO employs a nonparametric KDE approach, which offers greater modeling flexibility. However, it exhibits unstable estimates of $z(x)$ and wide confidence intervals, rendering it less reliable than its parametric counterparts.  

RSOPSO accounts for regime switching but suffers from parameter volatility. In the 4-regime exponential case, it often underestimates one rate while others remain large, leading to the poorest performance. Though more stable in the simpler second case, it remains unreliable and prone to large variations due to ignorance of input uncertainty.

RSOBSO captures a broad posterior range and effectively models multiple regimes, especially extreme ones like Regimes 1 and 4 in the first case. It consistently adapts well to regime switching. For instance, achieving lower GAPs during Regime 1 (Stages 114$\sim$119) and the abrupt shift to Regime 4 (Stages 119$\sim$124) in the first case, and smoothly adjusting through transitions from Regime 1 to 3 (Stages 53$\sim$55) and 2 to 1 (Stages 62$\sim$66) in the second. Overall, it demonstrates superior performance compared to all benchmark methods, with robust parameter inference and strong adaptability to dynamic regime shifts.

\subsection{Real-World Regime-Driven Inventory Problems}
\label{Sec: inventory_problem}
The 2008-2009 Great Recession presented substantial challenges for inventory management, primarily driven by significant variability in consumer demand \citep{dooley2010inventory}. To further demonstrate the practical robustness of RSOBSO, we apply it to a real-world inventory problem characterized by regime-driven demands during this turbulent period. In contrast to synthetic experiments with controlled regime structures and known transitions, this application involves market-induced regime shifts with unknown transition probabilities, thereby providing a more realistic and stringent evaluation of the algorithm effectiveness.

Our case study is built on a classical $(s, S)$ inventory model, where a company periodically manages stock with a reorder point $s$ and an order-up-to level $S$. Orders are only placed when inventory falls below $s$, replenishing stock up to $S$ \citep{fu1997techniques}. In the regime-driven setting, each regime corresponds to a distinct demand distribution, necessitating the algorithm accurately identify the prevailing regime and adapt decisions in a responsive and informed manner.

To simulate realistic market dynamics, we adopt regime classifications from \citet{pun2023data}, which include a two-regime bull-bear market model and a four-regime framework developed by Bridgewater Associates based on a $2 \times 2$ interaction between U.S. CPI and GDP trends.

\subsubsection{2-Regime Inventory Problem}
\label{Sec: 2-inventory}
Customer demand is higher in bull market (Regime 1) and lower in bear market (Regime 2), modeled by exponential distributions with rates $\lambda_1 = 1/20$ and $\lambda_2 = 1$. Initial regime data was generated to span Jan 2004$\sim$Dec 2007 and additional data cover Jan 2008$\sim$Dec 2009 (Great Recession). 

\subsubsection{4-Regime Inventory Problem}
\begin{figure}[h]
	\centering
	\subfloat{\includegraphics[width=0.46\textwidth]{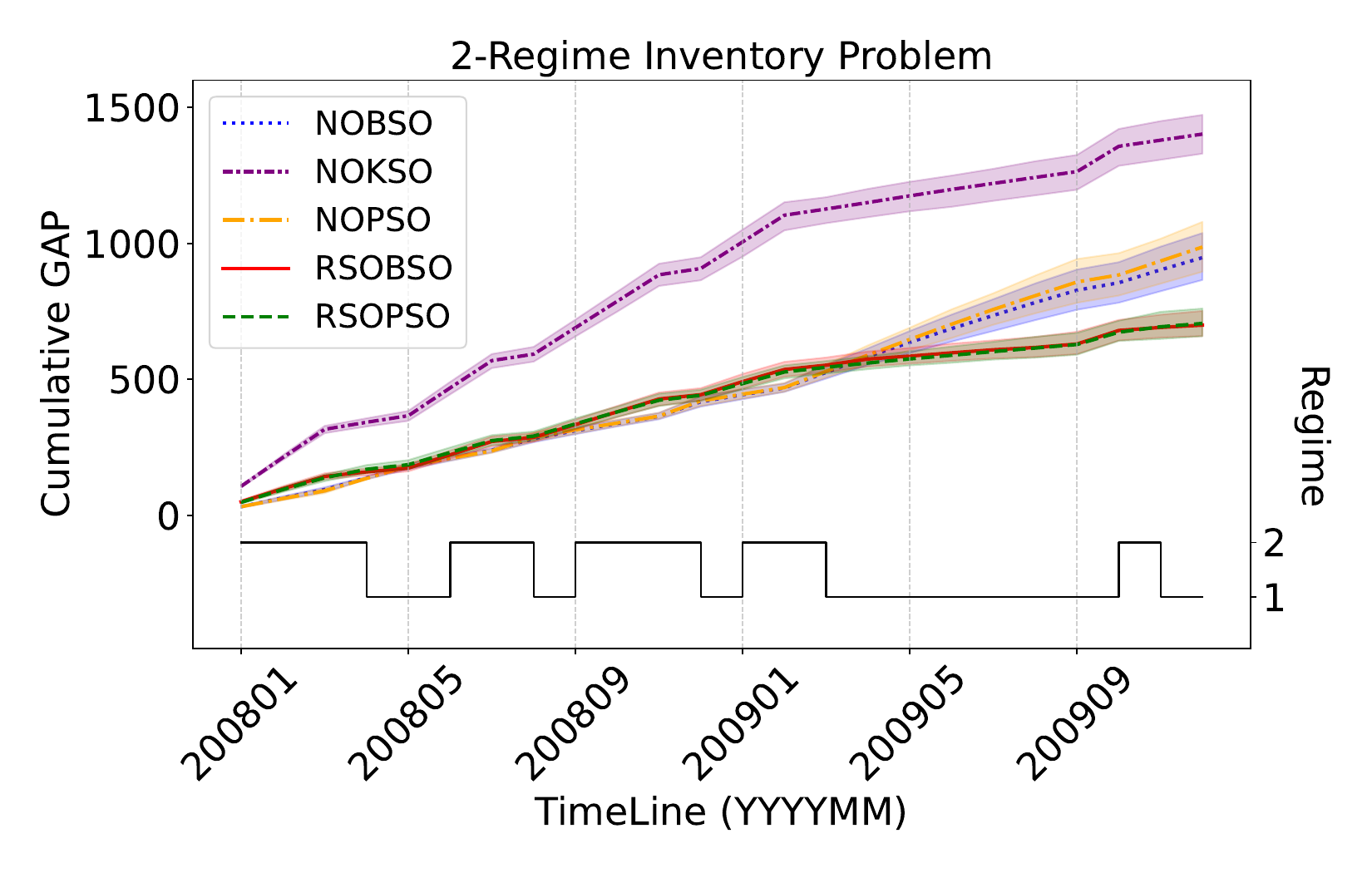}}
	\hfill
	\subfloat{\includegraphics[width=0.46\textwidth]{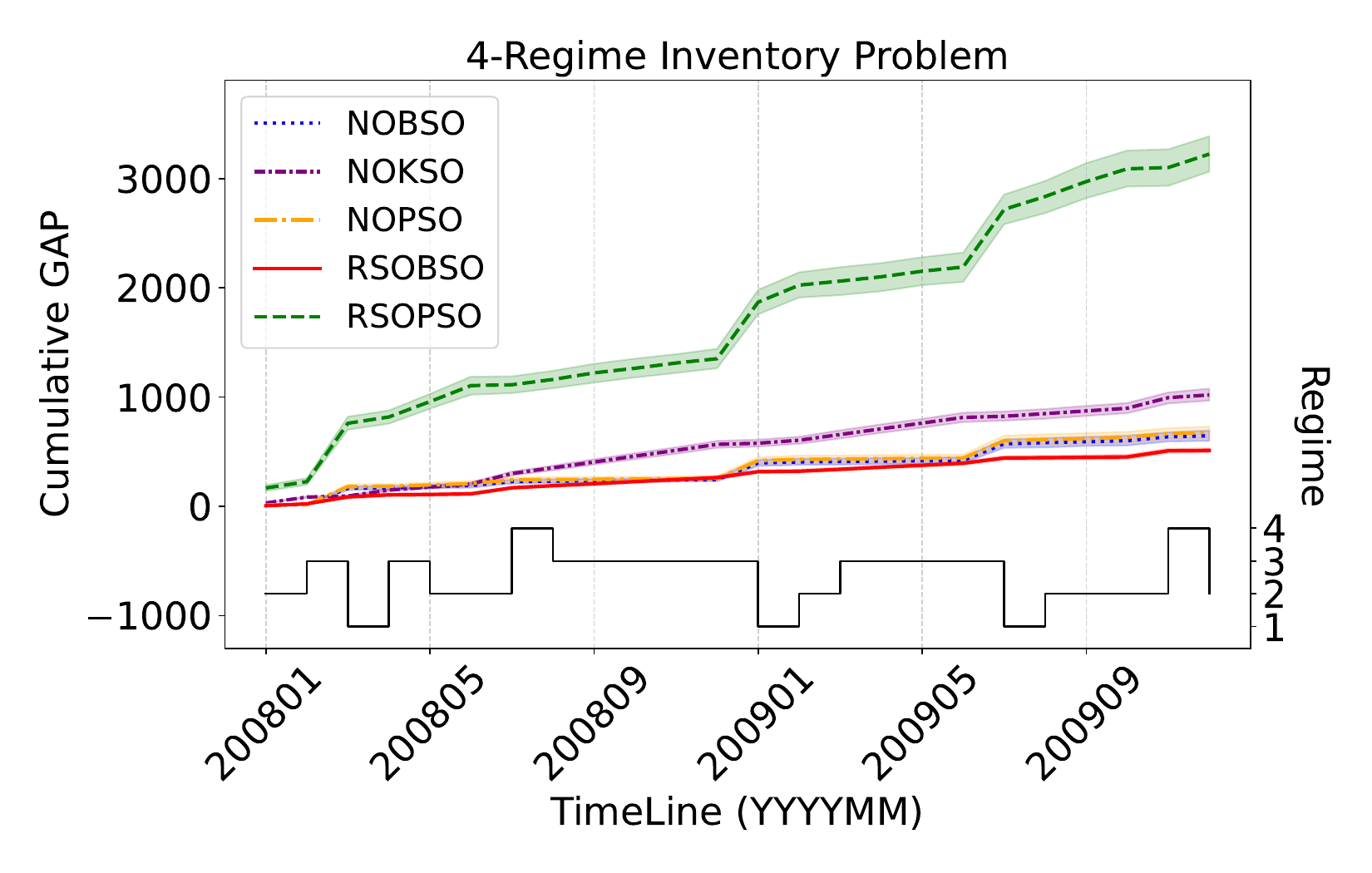}}
	\caption{Cumulative GAP for regime-driven inventory problems (black line indicates regimes).}
	\label{fig:combined_gap_inventory}
\end{figure}

As noted in \citet{mankiw2019macroeconomics}, high GDP growth boosts income and consumption, while high inflation undermines real income and dampens demand. Accordingly, consumer demand is expected to decrease across the following four regimes: high GDP \& low CPI, high GDP \& high CPI, low GDP \& low CPI, and low GDP \& high CPI. The corresponding exponential rates are $\{1/30, 1/18, 1/12, 1\}$. Initial regime data span Jan 2000$\sim$Dec 2007 and additional data cover Jan 2008$\sim$Dec 2009 (Great Recession). We plot the real regimes and the cumulative GAP of the two cases in Figure \ref{fig:combined_gap_inventory}.

\subsubsection{Analysis of the Regime-Driven Inventory Experiment}
The analysis of non-regime-switching methods aligns with findings from the synthetic problems. Consequently, we focus on regime-switching methods to assess how well they capture dynamic regime behaviors and underlying economic mechanisms during the 2008$\sim$2009 Great Recession.

In the two-regime setting, RS-based approaches successfully recover well-separated regime parameters, thereby achieving smaller cumulative GAPs and narrowing the performance gap between RSOBSO and RSOPSO. Prior to February 2009, frequent regime transitions hinder the predictive accuracy of RS-based methods. However, from that point onward, the system predominantly resides in Regime 1, during which RS-based approaches consistently yield the lowest GAP values. From an economic perspective, Regime 1 corresponds to a bull market, and its sustained prevalence post-2009 may indicate early signs of recovery from the Great Recession.

In the four-regime setting, RSOPSO has at least three large plug-in rate parameters and result in the weakest performance given the inventory objective heightened sensitivity to small rate values. This outcome is consistent with the first synthetic problem. RSOBSO, in contrast, accurately estimates three small and one large rate parameter, closely aligning with the underlying regime structure. This leads to markedly lower GAP in Regime 1, improved performance in Regimes 2 and 3, and competitive results in Regime 4. Although Regime 1 appears only three time, coinciding with key turning points also captured by non-RS methods, the dominance of Regimes 2 and 3 over a 19-month span accounts for the relatively modest overall GAP difference. From an economic standpoint, the infrequent occurrence of Regime 1, characterized by high GDP and low CPI, is consistent with its rarity during periods of economic recession.

In summary, RSOBSO emerges as the most robust method, effectively capturing regime-dependent uncertainty and delivering stable performance across all regimes. Notably, it achieves lower GAP even during the Great Recession, further highlighting its practical utility under volatile and extreme economic conditions.

\subsection{Real-World Portfolio Optimization Problems}
\label{Sec: portfolio}
To further illustrate the practical utility of RSOBSO, we apply it to a real-world portfolio optimization problem using two-asset return data from the 2008$\sim$2009 Great Recession. In contrast to the inventory setting, which assumes fully specified regime structures and known input distributions, this application relies solely on observed stochastic returns. The presence of noisy and high-dimensional inputs significantly increases the complexity of the problem, challenging the model's capacity to uncover underlying structure and make optimal decisions.  

\begin{figure}[h]
	\centering
	\subfloat{\includegraphics[width=0.46\textwidth]{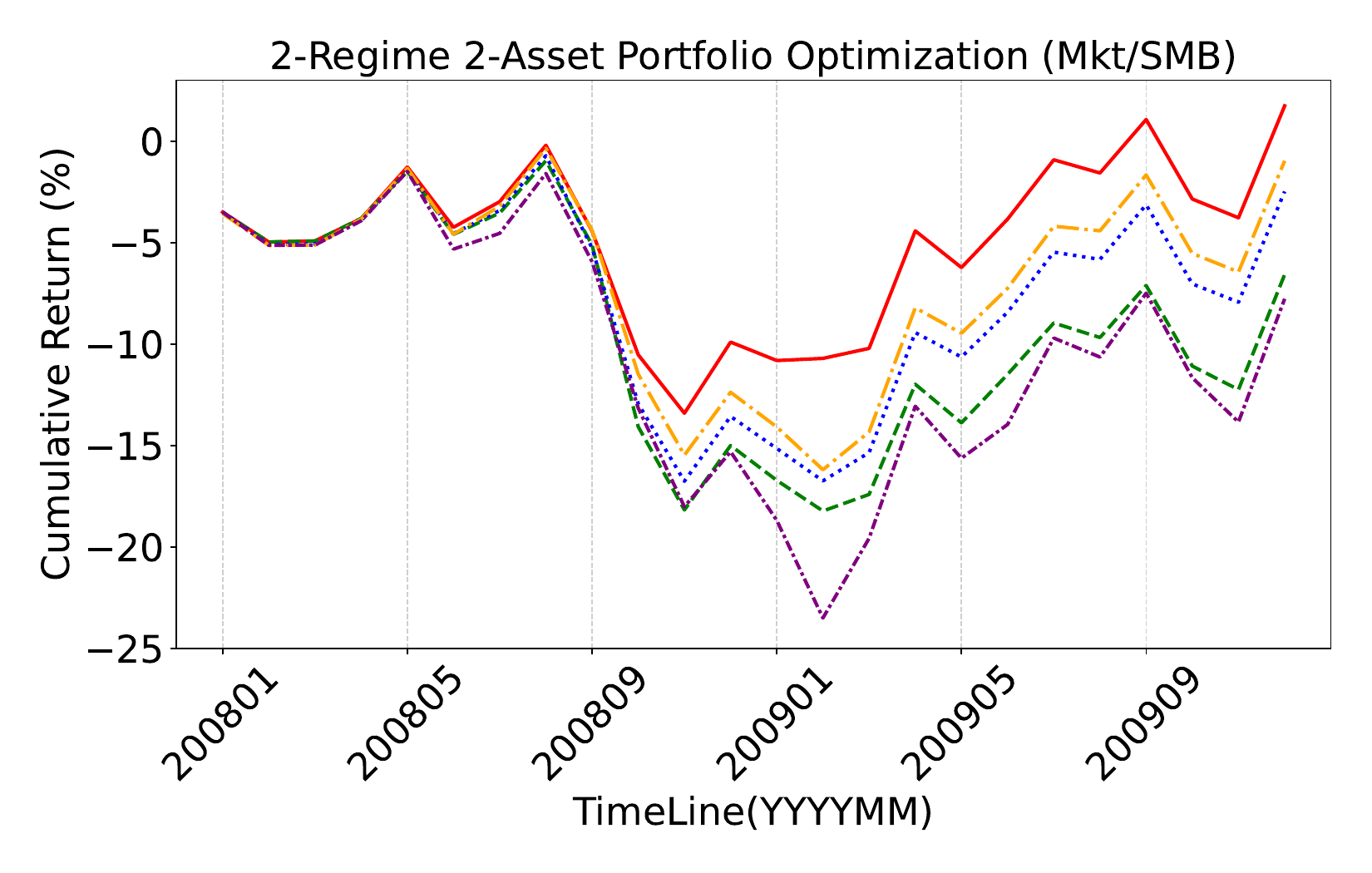}}
	\hfill
	\subfloat{\includegraphics[width=0.46\textwidth]{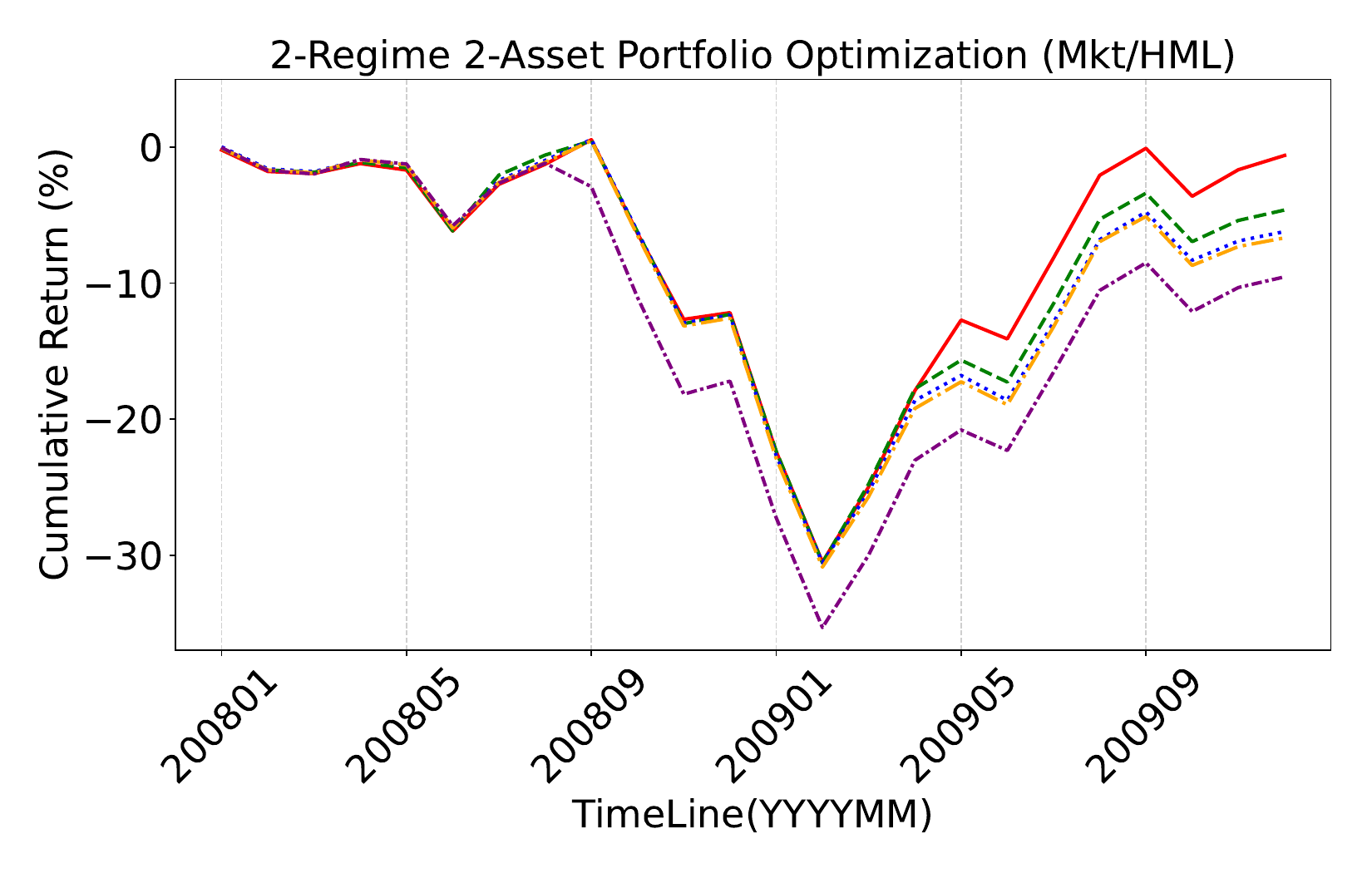}}
	\caption{Cumulative returns for 2-asset portfolio problems.}
	\label{fig:combined_gap_portfolio}
\end{figure}

The decision vector $x = (x_1, x_2)$ represents the portfolio weight, subject to $x_1, x_2\geq 0$ and $x_1+x_2 = 1$. Given the random return vector $r = (r_1, r_2)$, the resulting portfolio return is computed as $r^* = x\cdot r^T$. For a fixed decision $x$, the objective is the certainty equivalent (CEQ) return, which represents the risk-free return that an investor would accept in place of taking on a particular risky portfolio \citep{pun2023data}. Since the true return distribution is unknown, we approximate the CEQ 
using sample estimates of the mean and standard deviation of $r^*$, denoted by $\hat{\mu}(r^*)$ and $\hat{\sigma}(r^*)$, respectively. The estimated CEQ is then computed as:
$\widehat{\text{CEQ}}(x) = \hat{\mu}(r^*) - \frac{1}{2} \hat{\sigma}^2(r^*).$
We use the monthly return data of the Mkt/SMB/HML portfolios, as described in \citet{demiguel2009optimal}, with the initial sample period spanning from January 2004 to December 2007. For simplicity, each experiment considers only two assets. The emission distribution is assumed to be a bivariate Gaussian with a diagonal covariance matrix, corresponding to independent components with unknown means and variances. The cumulative return for the Mkt/SMB and Mkt/HML portfolios are presented in Figure~\ref{fig:combined_gap_portfolio}. 

\subsubsection{Analysis of Portfolio Optimization Experiments}
Under both the Mkt/SMB and Mkt/HML settings, RSOBSO consistently delivers the highest cumulative returns throughout the entire evaluation horizon. In the Mkt/SMB case, it demonstrates strong downside protection by allocating a substantial portion of the portfolio to SMB during the turbulent months from October 2008 to February 2009, precisely when market volatility peaks and other benchmark methods exhibit greater losses. This allocation effectively reduces exposure to market risk, resulting in a smoother return trajectory during the crisis. In the Mkt/HML case, RSOBSO shows high sensitivity to structural changes by quickly identifying the market recovery in May 2009 and shifting the allocation towards the Mkt asset to capitalize on the rebound. This timely adjustment enables it to outperform other methods that react more slowly or maintain static portfolios.

While most benchmark methods exhibit broadly similar portfolio trends, their delayed responses to abrupt return fluctuations can hinder timely decision-making, leading to missed opportunities and inadequate risk mitigation. In contrast, RSOBSO's ability to rapidly identify structural shifts in asset returns and update forecasts in real time enables more agile and informed portfolio adjustments, that said, an essential capability for navigating volatile financial market conditions and safeguarding investment performance, making it a powerful tool for real-world investment decisions.

\section{Extensions to Unknown Number of Regimes}
\label{sec:Extension}
This section presents an extension of RSOBSO to scenarios where the number of regimes in MSM is unknown. To address this, we employ the Hierarchical Dirichlet Process HMM (HDP-HMM) that offers a data-driven Bayesian nonparametric approach by placing an HDP prior on the regime transition distribution, and propose an algorithm (Algorithm 2 of EC.1) for inference of the number of regime. The details of the Bayesian inference are provided in EC.1. By combining Algorithm~\ref{algo-BRO} and 2 in an online setting, we develop the HDP-HMM RSOBSO algorithm (Algorithm 3), as detailed in EC.2. The algorithm dynamically decides an upper limit, $N_{\max}$, of regimes. Upon receiving a new observation, $N_{\max}$ is updated to $S_{\max}+1$ ($S_{\max}$ is the maximum inferred number of regimes in last stage) to account for the potential emergence of a previously unseen regime. 

\begin{figure}[h]
	\centering
	\includegraphics[width=0.5\textwidth]{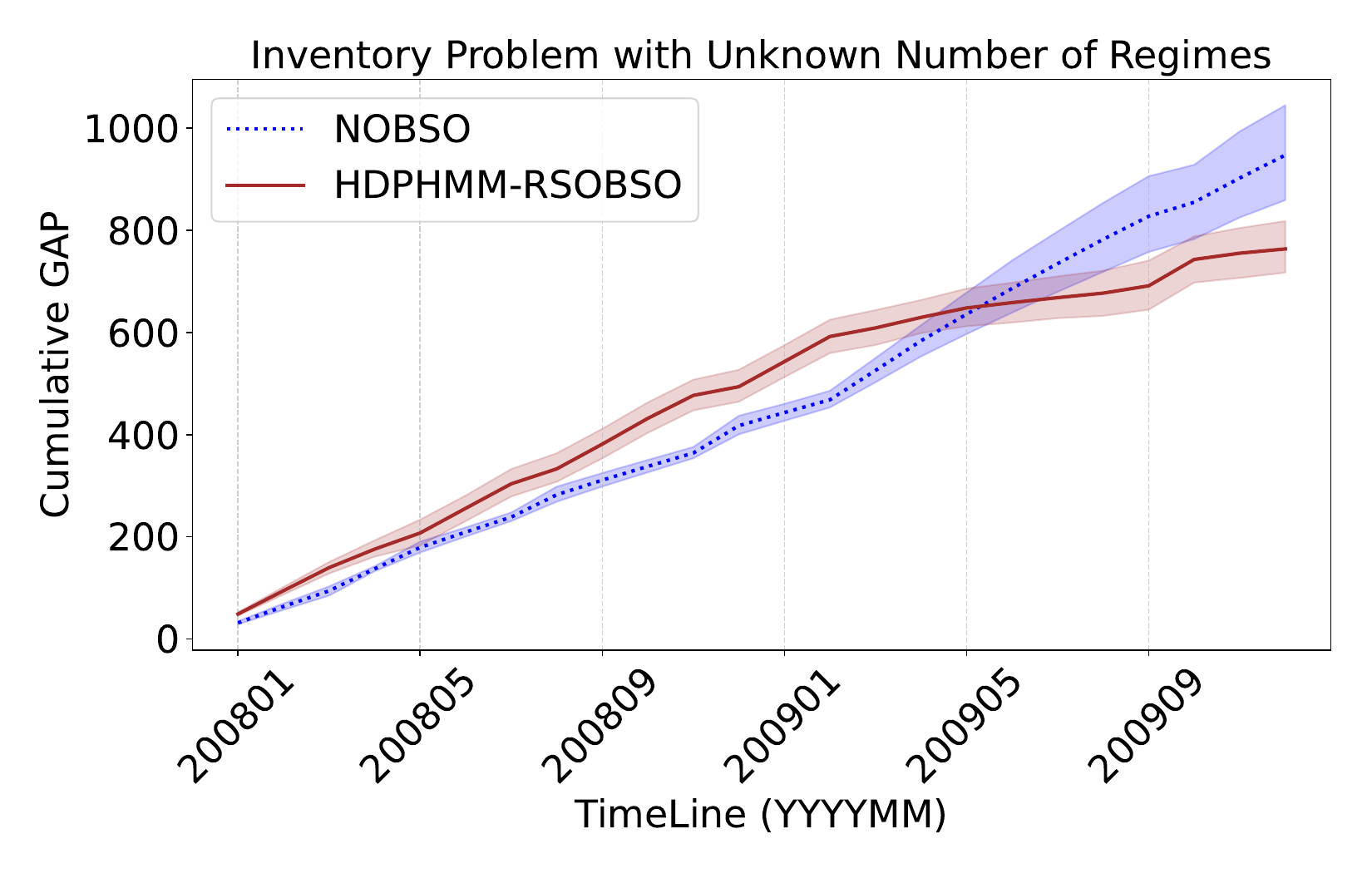}
	\caption{Cumulative GAP for extended 2-regime inventory problem.}
	\label{fig:combined_gap_states}
\end{figure}

Figure~\ref{fig:combined_gap_states} compares HDP-HMM RSOBSO with NOBSO in the 2-regime inventory experiment (Section~\ref{Sec: 2-inventory}). 
The inferred number of regimes falls between 2 and 4, closely aligning with the true number. For the 3-regime case, the estimated input parameters typically include two small values ($<0.1$) and one large value ($\sim$1), preserving a clear separation in magnitude. In the 4-regime case, two representative patterns emerge: (i) two small and two large values, resembling the 3-regime structure with one large component further split; (ii) three small and one large value, resembling the 2-regime case with the small component further divided. Despite minor deviations, parameter estimates remain distinctly separated in magnitude. Overall, HDP-HMM RSOBSO consistently outperforms NOBSO, showing superior robustness and efficiency under regime uncertainty.

\section{Summary}\label{sec:Conclusion}
This work investigates online simulation optimization with non-stationary streaming input data characterized by regime-switching distributions. We employ a Markov Switching Model (MSM) to capture non-stationary dynamics of the input data. A Bayesian framework is proposed to approximate the unknown true objective function, accounting for both input parameter uncertainty and the uncertainty in predicting regime-switching dynamics. We provide theoretical validation by establishing the consistency and asymptotic normality of the approximated objective function. To optimize this function, we integrate a joint Gaussian Process (GP) model to aggregate simulation results across sequential decision stages and utilize a regime-aware expected improvement function to select design points and input distributions for experimentation. Numerical experiments across various scenarios demonstrate that the proposed RSOBSO algorithm outperforms benchmark methods that overlook non-stationarity, input uncertainty, or joint GP modeling. Moreover, RSOBSO exhibits robust performance even when the number of regimes is unknown, highlighting its adaptability and effectiveness.


\appendix
 \section{Algorithm for Hidden Regime Number Inference}\label{EC: algo-HDPHMM}

 The stick-breaking construction of HDP-HMM is given below \citep{Teh01122006}:
 \[
 \begin{aligned}
	 \beta &\sim \text{GEM}(\gamma),\quad \pi_j \sim \text{DP}(\alpha, \beta),\quad S_t \sim \pi_{S_{t-1}},\\
	 \lambda_j &\sim H,\quad \xi_t \sim \tilde{P}(\lambda_{S_t}),\quad j = 1, 2, \ldots,\quad t = 1, 2, \ldots, T.
	 \end{aligned}
 \]
 Here, $\beta$ is the global probability vector over a potentially infinite number of hidden regimes, drawn from a stick-breaking process known as the GEM distribution with concentration parameter $\gamma$. Each regime-specific transition distribution $\pi_j$, the $j$-th row of the transition matrix, is drawn from a DP with concentration parameter $\alpha$ and base measure $\beta$, allowing regime transitions to share statistical strength via the shared global weights. 
 The emission parameter $\lambda_j$ for regime $j$ is drawn from the prior distribution $H$. 

 To conduct inference under this hierarchical model, we adopt the Blocked Gibbs sampler with a weak limit approximation to the DP transition prior: $\beta \sim \text{Dir}(\gamma/L,\dots,\gamma/L)$ and $\pi_j \sim \text{Dir}(\alpha\beta_1,\dots,\alpha\beta_L)$, where $L$ is set to exceed the expected number of regimes \citep{fox2011sticky}. This inference procedure is adapted to our setting, as described in Algorithm~\ref{algo-infer-states}. To estimate the number of regimes, we apply a threshold to the stick-breaking weights and compute the mode of the sampled regime counts after resampling.

 \begin{algorithm}
	 \caption{Inference of Hidden Regime Number via HDP-HMM}
	 \label{algo-infer-states}
	 \begin{algorithmic}[1]
		 \State \textbf{Input:} maximum number of regimes $N_{\max}$, sampling steps $M$, observations $\xi$, threshold $\tau$;
		 \State Initialize empty list $\mathcal{S} \leftarrow []$; 
		 \For{$i = 1$ to $M$}
		     \State Resample model parameters using WeakLimit HDP-HMM with $N_{\max}$ components and $\xi$; 
		     \State Extract stick-breaking weights $\beta$;
		     \State Count $s_i = \sum_{j=1}^{N} \mathbb{I}(\beta_j \geq \tau)
		 $ and append $s_i$ to $\mathcal{S}$;
		 \EndFor
		 \State $S_{\max} \leftarrow \max(\mathcal{S})$, $S_{\text{mode}} \leftarrow$ mode of $\mathcal{S}$;
		 \State \textbf{Output:} maximum inferred number of regimes $S_{\max}$, and inferred number of regimes $\hat{R} \leftarrow S_{\text{mode}}$;
		 \end{algorithmic}
	 \end{algorithm}

 \section{HDPHMM-RSOBSO Algorithm}\label{EC: HDPHMMRS}
 The proposed HDPHMM-RSOBSO Algorithm is shown in Algorithm \ref{HDPHMM-RS-BRO}.
 \begin{algorithm}
	 \caption{The HDPHMM-RSOBSO Algorithm}
	 \label{HDPHMM-RS-BRO}
	 \begin{algorithmic}[1]
		 \State Given $N_{max}$, $M$, threshold $\tau_0$, $\xi^0$ of size $h$, calculate $S_{max}$ and $\hat{R}_0$ using Algorithm \ref{algo-infer-states};  
		 \State Generate $\{x_1,\cdots,x_{n_0}\}$ using LHS; given prior $P_0(\vartheta)$, draw $n_0$ sets of posterior samples from $P(\vartheta|\xi^0)$ to obtain initial points $\{(x_1, {\lambda_1^{(1)}}),\cdots,(x_{n_0}, {\lambda_1^{(n_0)}}), \cdots, (x_1, {\lambda_{\hat{R}_0}^{(1)}}),\cdots,(x_{n_0}, {\lambda_{\hat{R}_0}^{(n_0)}})\}$; 
		 \State Set $n \leftarrow n_0 * \hat{R}_0$;
		 \State Run the simulation experiments at initial points with $m$ replications each and get observed output sample mean vector ${\bar{Y}^{(0)}_n}$; 
		 \State \textbf{GP fitting:} Based on $\bar{Y}^{(0)}_n$, construct GP model $Z^{(0)}_n({x}, \lambda)$, generate $N_{MC}$ sets of samples from $P(\vartheta|\xi^0)$, and derive $\hat{G}_0({x})$;
		 \For{$t = 0$ to $t_{max}-1$}
		 \For{$e_t = 0$ to specified number of optimization steps $u-1$}
		 \State Select the next evaluation point $({x}^{(t)}_{n+e_t+1}, \lambda^{(t,n+e_t+1)})$ using the refined EGO algorithm;
		 \State Run simulation experiments at $({x}^{(t)}_{n+e_t+1}, \lambda^{(t, n+e_t+1)})$ with $m$ replications and obtain  $\bar{y}({x}^{(t)}_{n+e_t+1}, \lambda^{(t, n+e_t+1)})$, \text{\hspace{2.8em}} set ${\bar{Y}^{(t)}_{n}} \leftarrow [ {\bar{Y}_{n}^{(t)}},\bar{y}({x}^{(t)}_{n+e_t+1}, \lambda^{(t, n+e_t+1)})]^T$; 
		 \State Set $n \leftarrow n+e_t+1$;
		 \State \textbf{Update}: update GP model based on ${\bar{Y}^{(t)}_{n}}$; 
		 \EndFor
		 \State \textbf{return} $\hat{x}_t^* \leftarrow \arg\min_{x}\mu_t({x})$;
		 \State Given streaming input data $\xi^{t+1} \leftarrow \xi^t \cup \xi_{t+1}$, \textbf{update} $N_{max} \leftarrow S_{max} + 1$; 
		 \State Calculate $S_{max}$ and $\hat{R}_{t+1}$ with threshold $\tau_{t+1}$ updated $N_{max}$ and $\xi^{t+1}$ through Algorithm \ref{algo-infer-states};  
		 \State Generate $N_{MC}$ samples from $P(\vartheta|\xi^{t+1})$, use the latest updated $Z_n^{(t)}({{x}}, \lambda)$ to derive $\hat{G}_{t+1}({x})$;   
		 \State Set $\bar{Y}_{n}^{(t+1)} \leftarrow \bar{Y}_{n}^{(t)}$.
		 \EndFor
		 \end{algorithmic}
	 \end{algorithm}

 \section{Proofs for Section \ref{sec:Theory}}
 \label{EC: theory_proof}
 \begin{proof}{Proof of Lemma \ref{cons_msm}}
	 Assumptions (B1)$\sim$(B6) in \citet{de2008asymptotic} ensure the strong consistency of MLE over any compact parameter space for MSM with a finite regime space. In this context, Assumptions (B1)-(B2) in \citet{douc2020posterior} are satisfied with the parameter space $\Theta$. Meanwhile, Assumption \ref{assum2} guarantees that Assumptions (B3)-(B4) in \citet{douc2020posterior} hold. Therefore, all the assumptions for Theorem 1 in \citet{douc2020posterior} are satisfied and the lemma follows.   
	\end{proof}

 \begin{proof}{Proof of Theorem \ref{consistency_obj}}
	 From the derivations of MSM, we have:
	 $$P(S_{t+1}=l|\xi^t,\vartheta)=\sum_{k=1}^{\tilde{R}} A_{k,l}P(S_t=k|\xi^t,\vartheta),$$
	 where $$P(S_t=k|\xi^t,\vartheta)=\frac{\tilde{P}(\xi_t|\lambda_k)P(S_t=k|\xi^{t-1},\vartheta)}{\sum_{j=1}^{\tilde{R}} \tilde{P}(\xi_t|\lambda_j)P(S_t=j|\xi^{t-1},\vartheta)}.$$
	 The derivations are recursive and the denominator of $P(S_t=k|\xi^t,\vartheta)$ is not 0. By Assumption \ref{assum1} $(ii)$, it can be readily verified that $P(S_t=l|\xi^t,\vartheta)$ and $P(S_{t+1}=l|\xi^t,\vartheta)$ are continuous on $\Theta$ for $l=1,\cdots,\tilde{R}$. Moreover, we note that 
	 $$\widetilde{K}(x,\vartheta) =\mathbb{E}_{P(S_{t+1}|\xi^{t},\vartheta)}[\mathbb{E}_{P(\xi_{t+1}|S_{t+1},\vartheta)}(y(x,\xi_{t+1}))] =\sum_{i=1}^{\tilde{R}} P(S_{t+1}=i|\xi^{t},\vartheta)\mathbb{E}_{\tilde{P}(\xi_{t+1}|\lambda_{i})}[y(x,\xi_{t+1})].$$ 
	 When $x$ is fixed, we write $\widetilde{K}(x,\vartheta)$ simply as $\widetilde{K}(\vartheta)$. Given $S_{t+1}=i$, by Assumption \ref{assum1} $(i)$, 
	 $$\mathbb{E}_{\tilde{P}(\xi_{t+1}|\vartheta)}[y(x,\xi_{t+1})]=\mathbb{E}_{\tilde{P}(\xi_{t+1}|\lambda_{i})}[y(x,\xi_{t+1})]$$ 
	 is continuous on $\Theta$ for every fixed $x$. Thus, $\widetilde{K}$ is continuous on $\Theta$. Assumption (B1) in \citet{de2008asymptotic} gives the compactness of $\Theta$. Therefore, $\widetilde{K}$ is continuous and bounded on $\Theta$. The proof is completed by following from Lemma \ref{consistency_obj} and Definition 3.2 (weak convergence) in \citet{wu2018bayesian}. 
	 \end{proof}

 \begin{proof}{Proof of Theorem \ref{theorem2}}
	 Similar to the proof of Theorem \ref{consistency_obj}, it can be shown that $\mathbb{E}_{P_t}[\widetilde{K}(\cdot,\vartheta)]\rightarrow \widetilde{K}(\cdot, \vartheta^c)$ pointwise on $\mathcal{X}$ a.s. $(\eta)$ as $t\rightarrow\infty$. We only need to show $\mathbb{E}_{P_t}[\widetilde{K}(\cdot,\vartheta)]$ has a common Lipschitz constant $L$ for all $t$. Then from Lemma 3.7 in \citet{wu2018bayesian}, $\mathbb{E}_{P_t}[\widetilde{K}(\cdot,\vartheta)]\rightarrow \widetilde{K}(\cdot,\vartheta^c)$ uniformly on $\mathcal{X}$, and the proof is completed by following Lemma 3.6 in \citet{wu2018bayesian}. Let $\mathcal{F}_t := \sigma\{(\xi_k), k\leq t\}$ denote the filtration generated by data and $\mathcal{F}_{\infty} := \sigma(\bigcup\limits_t \mathcal{F}_t)$. For $\forall x,y \in \mathcal{X}$,
	 \begin{equation*}
		 \begin{aligned}
			 \left| \mathbb{E}_{P_t}[\widetilde{K}(x,\vartheta)] - \mathbb{E}_{P_t}[\widetilde{K}(y,\vartheta)] \right|
			 &\leq \mathbb{E}_{P_t}\left| \widetilde{K}(x,\vartheta) - \widetilde{K}(y,\vartheta) \right| \\
			 &\leq \mathbb{E}_{P_t}[\psi(\vartheta) \Vert x - y \Vert] \\
			 &= \mathbb{E}_{P_t}[\psi(\vartheta)] \Vert x - y \Vert \\
			 &= \mathbb{E}_{\eta}[\psi(\vartheta) \mid \mathcal{F}_t] \Vert x - y \Vert.
			 \end{aligned}
		 \end{equation*}
	 Since $$\mathbb{E}_{\eta}\left|\psi(\vartheta)\right|=\mathbb{E}_{\eta}[\psi(\vartheta)]=\int_{\Theta} \psi(\vartheta)P_0(d\vartheta)<\infty,$$
	 then $\mathbb{E}_{\eta}[\psi(\vartheta)|\mathcal{F}_t]$ is a Doob martingale. By the martingale convergence theorem (Theorem 4.6.8 of \citep{durrett2019probability}), $\mathbb{E}_{\eta}[\psi(\vartheta)|\mathcal{F}_t]\rightarrow \mathbb{E}_{\eta}[\psi(\vartheta)|\mathcal{F}_{\infty}]$ as $t\rightarrow \infty.$ Since 
	 $$\mathbb{E}_{\eta}\{\mathbb{E}_{\eta}[\psi(\vartheta)|\mathcal{F}_{\infty}]\}=\mathbb{E}_{\eta}[\psi(\vartheta)]<\infty,$$ 
	 $\mathbb{E}_{\eta}[\psi(\vartheta)|\mathcal{F}_\infty]$ is a.s. finite. Thus, there exists a Lipschitz constant $L := \sup_t \mathbb{E}_{\eta}[\psi(\vartheta)|\mathcal{F}_t]<\infty, $ and the proof is completed.  
	 \end{proof}

 \begin{proof}{Proof of Lemma \ref{lemma-added}}
	 Assumptions (B1)-(B6) and (A5) from \citep{de2008asymptotic} imply the Central Limit Theorem for the score function at $\vartheta^c$: $$\frac{1}{\sqrt{t}}\nabla L_t(\vartheta^c) \Rightarrow \mathcal{N}(0, J_c),$$ as $t\rightarrow\infty$. 
	 Thus, we have $\Delta_t \Rightarrow \mathcal{N}(0,J_c^{-1})$ as $t\rightarrow\infty$. Let $p_{Z_t}$ be the density of $P_{Z_t}$, and $\phi^{**}$ be the Gaussian density with mean $\Delta_t$ and covariance matrix $J_c^{-1}$, $\Phi^{**}$ be the corresponding Gaussian distribution. We have $\int_\Theta |p_{Z_t}(s)-\phi^{**}(s)|ds \rightarrow0.$ Then, for $A\in \Theta$,
	 \begin{equation*}
		 \begin{aligned}
			      |P_{Z_t}(A)-\Phi^{**}(A)|
			      &=\left|\int_A p_{Z_t}(s)\,ds - \int_A \phi^{**}(s)\,ds\right|\\
			      &\leq \int_A \left|p_{Z_t}(s) - \phi^{**}(s) \right|\,ds\\
			      &\leq \int_\Theta \left|p_{Z_t}(s) - \phi^{**}(s) \right|\,ds.
			 \end{aligned}
		 \end{equation*}
	 Therefore, 
	 $$\sup_{A\in\Theta} |P_{Z_t}(A)-\Phi^{**}(A)|\leq \int_\Theta \left|p_{Z_t}(s)-\phi^{**}(s) \right|ds\rightarrow0.$$ 
	 Thus, $\Vert P_{Z_t}-\Phi^{**}\Vert_\textup{TV}\rightarrow 0$, equivalently, we have $\Vert P_{Z_t}-\mathcal{N}(\Delta_t, J_c^{-1})\Vert_\textup{TV}\rightarrow 0$ in probability $(P_{\vartheta^c}^N)$ as $t\rightarrow \infty$.   
	 \end{proof}

 \begin{proof}{Proof of Theorem \ref{mean-normality}}
	 The first-order Taylor expansion of $\widetilde{K}$ around $\vartheta^c$ yields 
	 $$\mathbb{E}_{P_t}[\sqrt{t}(\widetilde{K}(\vartheta)-\widetilde{K}(\vartheta^c))]=\nabla \widetilde{K}(\vartheta^c)^T \mathbb{E}_{P_t}[\sqrt{t}(\vartheta-\vartheta^c)]+\mathbb{E}_{P_t}[e(\vartheta) \Vert \sqrt{t}(\vartheta-\vartheta^c) \Vert],$$
	 where $e(\vartheta)\rightarrow 0$ if $\vartheta\rightarrow \vartheta^c$. Since the Bayes estimator is the mean of posterior distribution, Theorem 2.2 in \citet{de2008asymptotic} establishes that $\sqrt{t}(\mathbb{E}_{P_t}(\vartheta)-\vartheta^c)\Rightarrow \mathcal{N}(0,J_c^{-1})$. Therefore, the term $\nabla \widetilde{K}(\vartheta^c)^T \mathbb{E}_{P_t}[\sqrt{t}(\vartheta-\vartheta^c)]$ converges weakly to $\mathcal{N}(0,\sigma_x^2).$ For the term $\mathbb{E}_{P_t}[e(\vartheta) \Vert \sqrt{t}(\vartheta-\vartheta^c) \Vert],$ we apply the H\"older's inequality: 
	 $$\left| \mathbb{E}_{P_t}[e(\vartheta) \Vert \sqrt{t}(\vartheta-\vartheta^c) \Vert]\right|\leq (\mathbb{E}_{P_t}[|e(\vartheta)|^{\frac{1+\gamma}{\gamma}}])^{\frac{\gamma}{1+\gamma}} (\mathbb{E}_{P_t}[\Vert \sqrt{t} (\vartheta-\vartheta^c)\Vert ^{1+\gamma}])^{\frac{1}{1+\gamma}}.$$
	 Since $\widetilde{K}$ is continuous and $\Theta$ is compact, we assume $e(\cdot)$ is bounded and continuous on $\Theta$. And assuming $e(\vartheta^c)=0$ does not affect the above Taylor expansion terms. By Lemma \ref{cons_msm} and Definition 3.2 (weak convergence) in \citet{wu2018bayesian}, as $t\rightarrow \infty$, $$\mathbb{E}_{P_t}[|e(\vartheta)|^{\frac{1+\gamma}{\gamma}}]\rightarrow |e(\vartheta^c)|^{\frac{1+\gamma}{\gamma}}=0,$$ a.s. $(P_{\vartheta^c}^N)$.
	 From Assumption 4.1 in \citet{wu2018bayesian}, we obtain that $\mathbb{E}_{P_t}[\Vert \sqrt{t} (\vartheta-\vartheta^c)\Vert ^{1+\gamma}]$ is bounded in probability $(P_{\vartheta^c}^N)$. It implies that $\mathbb{E}_{P_t}[e(\vartheta) \Vert \sqrt{t}(\vartheta-\vartheta^c) \Vert]$ converges weakly to 0. Therefore, the proof is completed. 
	 \end{proof}

 \begin{proof}{Proof of Theorem \ref{opt-normality}}
	 The whole proof process is divided into three steps.
	
	 Step 1: By the Tychonoff product theorem in \citet{royden2010real}, $\mathcal{X}\times \Theta$ is compact. Then by the Heine-Cantor theorem (Theorem 4.19 in \citep{rudin1964principles}), $\widetilde{K}$ is uniformly continuous on $\mathcal{X}\times \Theta$. Thus, for $\forall \epsilon>0$, there exists $\delta>0$ such that $|\widetilde{K}(x,\vartheta)-\widetilde{K}(x^\prime, \vartheta^\prime)|<\epsilon$ as long as $\Vert (x,\vartheta)-(x^\prime, \vartheta^\prime) \Vert <\delta$. If $\Vert (x,\vartheta)-(x^\prime, \vartheta) \Vert = \Vert(x-x^\prime) \Vert <\delta$, by Lemma 3.8 in \citet{wu2018bayesian}, 
	 $$|\mathbb{E}_{P_t}\{\widetilde{K}(x,\vartheta)\}-\mathbb{E}_{P_t}\{\widetilde{K}(x^\prime,\vartheta)\}|\leq \sup_{\vartheta\in\Theta} |\widetilde{K}(x,\vartheta)-\widetilde{K}(x^\prime, \vartheta^\prime)|<\epsilon.$$ Therefore, $\mathbb{E}_{P_t}\{\widetilde{K}(\cdot,\vartheta)\}$ is uniformly continuous. 
	
	 Step 2: This step is to show the weak convergence of finite-dimensional distributions. Since the expectation is a linear functional, we have $$g_t(x) = \mathbb{E}_{P_t}\{\sqrt{t}[\widetilde{K}(x,\vartheta)-\widetilde{K}(x,\vartheta^c)]\}.$$ 
	 Fix a finite sequence $x_1, x_2, \cdots, x_k \in \mathcal{X}$. For $[g_t(x_1), g_t(x_2),\cdots, g_t(x_k)]$, we apply a Taylor expansion inside the expectation for each dimension. The remainder terms also form a $k$-dimensional random vector, which converges in probability to 0 if and only if each dimension does. Thus, the proof of Theorem \ref{mean-normality} can be easily extended to show each formulation's finite-dimensional distributions converge weakly to that of $Y_x$.
	
	 Step 3: By Theorem 7.5 in \citet{billingsley2013convergence}, the proof will be completed if $g_t(\cdot)$ is stochastic equicontinuity (s.e.). The condition of s.e. is
	 $$\lim_{\delta\rightarrow0} \limsup_{t\rightarrow\infty} P_{\vartheta^c}^N(\zeta(g_t,\delta)\geq \epsilon)=0 \quad \forall \epsilon>0,$$
	 where $\zeta(f,\delta)$ for $f \in C(\mathcal{X})$ is defined as 
	 $\zeta(f,\delta) := \sup_{x,x^\prime\in \mathcal{X}, \Vert x-x^\prime\Vert<\delta, }|f(x)-f(x^\prime)|.$
	 By Taylor expansion, 
	 $$g_t(x) = \nabla_{\vartheta} \widetilde{K}(x, \vartheta^c)^T \mathbb{E}_{P_t}[\sqrt{t}(\vartheta-\vartheta^c)]+\mathbb{E}_{P_t}[e(x,\vartheta)\Vert \sqrt{t}(\vartheta-\vartheta^c)\Vert].$$ By Assumption 4.1 in \citet{wu2018bayesian}, $\mathbb{E}_{P_t}[\sqrt{t}(\vartheta-\vartheta^c)]$ converges weakly, for any $\iota>0$, there exists $M_\iota>0$ such that 
	 $P_{\vartheta^c}^N(\Vert \mathbb{E}_{P_t}[\sqrt{t}(\vartheta-\vartheta^c)]\Vert < M_\iota)>1-\iota \; $ $\forall t.$
	 Since we assume $\nabla_{\vartheta}\widetilde{K}(\cdot, \vartheta^c)$ is continuous (and hence uniformly continuous) on $\mathcal{X}$, for any $\epsilon>0$, there exists $\delta_{\iota}>0$ such that 
	 $$\sup_{x,x^\prime\in \mathcal{X}, \Vert x-x^\prime\Vert<\delta, } \Vert \nabla_\vartheta \widetilde{K}(x, \vartheta^c)-\nabla_\vartheta \widetilde{K}(x^\prime, \vartheta^c)\Vert < \frac{\epsilon}{M_\iota}.$$
	 On the event $\{\Vert \mathbb{E}_{P_t}[\sqrt{t}(\vartheta-\vartheta^c)] \Vert <M_\iota\}$, we have
	 $$\zeta(\nabla_{\vartheta} \widetilde{K}(x, \vartheta^c)^T \mathbb{E}_{P_t}[\sqrt{t}(\vartheta-\vartheta^c)], \delta_\iota)<\frac{\epsilon}{M_\iota}M_\iota=\epsilon,$$ thus the first term of $g_t(x)$ has the s.e. property. For the second term of $g_t(x)$, we only need to prove 
	 $$\sup_{x\in\mathcal{X}} |\mathbb{E}_{P_t}[e(x,\vartheta)\Vert \sqrt{t}(\vartheta-\vartheta^c)\Vert]|\Rightarrow 0.$$
	 Since $$\sup_{x\in\mathcal{X}} |\mathbb{E}_{P_t}[e(x,\vartheta)\Vert \sqrt{t}(\vartheta-\vartheta^c)\Vert]| \leq \mathbb{E}_{P_t}\left[\sup_{x\in\mathcal{X}}|e(x,\vartheta)|\Vert \sqrt{t}(\vartheta-\vartheta^c)\Vert\right],$$ it suffices to show for $\vartheta \sim P_t$ that $\sup_{x\in\mathcal{X}}|e(x,\vartheta)|\Rightarrow 0$ a.s. $(P_{\vartheta^c}^N)$. The continuity of $e$ on $\mathcal{X}\times \Theta$ implies $\sup_{x\in\mathcal{X}}|e(x,\cdot)|$ is continuous on $\Theta$. Setting $e(x,\vartheta^c)=0$ for all $x\in\mathcal{X}$ does not affect Taylor expansion. Therefore, $$\sup_{x\in\mathcal{X}}|e(x,\vartheta)|\Rightarrow \sup_{x\in\mathcal{X}}|e(x,\vartheta^c)|=0$$ a.s. $(P_{\vartheta^c}^N),$
	 and the rest follows from proof of Theorem \ref{mean-normality}. 
	 \end{proof}

 \section{Experiment Settings in Section \ref{sec:Experiment} and \ref{sec:Extension}}
 Detailed settings for all numerical experiments are presented in this section. As mentioned in Section \ref{Sec: GP}, a symmetric Dirichlet prior $\text{Dir}(1, \cdots, 1)$ is imposed on each row of the transition matrix throughout all experiments. 

 \label{Sec: parameter-settings}
 \subsection{4-Regime Exponential Emissions Case}
 The candidate space is $x \in [0, 50]$. The conjugate weakly informative prior distribution for $\lambda$ is $\mathrm{Gamma}(1, 0.1)$. The transition matrix $A$ is
 $$\begin{bmatrix}
	 0.7 & 0.1 & 0.1 & 0.1 \\
	 0.1 & 0.7 & 0.1 & 0.1 \\
	 0.1 & 0.1 & 0.7 & 0.1 \\
	 0.05 & 0.05 & 0.1 & 0.8
	 \end{bmatrix},
 $$
 $h = 100$, $u = 30$, $t_{max} = 25$, $m=100$, $N_{MC} = 100$.
 \subsection{3-Regime Gaussian Emissions Case}
 The candidate space is $\{x=(x_1, x_2)\mid x_1\in[-20, 20], x_2\in[-40, 40]\}$. The prior distribution for unknown $\mu$ is the uniform distribution $\mathcal{U}(0, 50)$. The transition matrix $A$ is
 $$\begin{bmatrix}
	 0.7 & 0.15 & 0.15 \\
	 0.15 & 0.7 & 0.15 \\
	 0.1 & 0.1 & 0.8 
	 \end{bmatrix},
 $$
 $h = 50$, $u = 30$, $t_{max} = 25$, $m=100$, $N_{MC} = 100$.
 \subsection{2-Regime Inventory Problem}
 In the inventory problem simulator, we set fixed ordering cost = 100, unit cost = 1, holding cost = 1, and back-order cost = 100 to compute the cost \citep{jalali2017comparison}. The candidate decision space is $\{x=[s,S]|s \in [1,69], S \in [70,250]\}$. The optimal solutions are $(63.8,127.0)$ with cost 147 for Regime 1, and $(1,70)$ with cost 38 for Regime 2. The prior distribution for unknown $\lambda$ is $\mathrm{Gamma}(1, 1)$. In this case, $h = 48$, $u = 30$, $t_{max} = 24$, $m=10$, $N_{MC} = 100$.
 \subsection{4-Regime Inventory Problem}
 We use the same values in the inventory problem simulator and candidate space as as 2-regime case. The optimal solutions and costs for the four regimes are as follows. Regime 1: (69, 191), cost = 222; Regime 2: (57, 118), cost = 135; Regime 3: (35, 87), cost = 97; Regime 4: (1, 70), cost = 38.The prior distribution for unknown $\lambda$ is $\mathrm{Gamma}(1, 0.1)$. In this case, $h = 96$, $u = 30$, $t_{max} = 24$, $m=10$, $N_{MC} = 100$. 

 \subsection{Portfolio Optimization Problems}
 Uniform priors are assigned to the unknown mean $\mu$ and standard deviation $\sigma$: $\mathcal{U}(0, 50)$ and $\mathcal{U}(0.1, 20)$ in the Mkt/SMB case, while $\mathcal{U}(0, 20)$ and $\mathcal{U}(0.1, 10)$ in the Mkt/HML case. In both cases, $h = 48$, $u = 30$, $t_{max} = 24$, $m=1000$, $N_{MC} = 100$.     

 \subsection{Unknown Number of Regimes}
 We set the prior $H = \mathrm{Gamma}(1, 1) $ for the exponential emission rates. The concentration parameters $\alpha$ and $\gamma$ follow $\mathrm{Gamma}(5, 1)$ priors. The threshold is set as $\tau_t = 1 / \sqrt{h+t}, t = 0, \cdots, t_{max}$. We use a truncation level $N_{\max} = 10$ and draw $M = 100$ posterior samples at each time stage.

 \section{Objective Functions of Benchmark Methods}\label{EC: obj}
 The benchmark methods adopt different forms of the objective function, as detailed below.
 \begin{itemize}
	 \item RSOPSO evaluates the objective using the weighted plug-in estimator: $\sum_{l =1}^{\tilde{R}}\hat{w}_l^t z(x, \hat{\lambda}_l^t)$, where $\hat{\lambda}_l^t$ is the plug-in posterior mean parameters and $\hat{w}_l^t$ is the posterior weight for regime $l$. 
	 \item NOBSO applies Bayesian averaging without regime switching: $\frac{1}{N_{MC}} \sum_{i = 1}^{N_{MC}}z(x, \lambda^{(t,i)})$, using MCMC samples $\lambda^{(t,i)}$ from the posterior. 
	 \item NOPSO uses a single plug-in estimator: $z(x,\hat{\lambda}^t)$, with $\hat{\lambda}^t$ being the posterior mean of the input parameter.
	 \item NOKSO is fully nonparametric and does not model input parameters explicitly, directly evaluating the objective $z(x)$.
	 \end{itemize}

 \section{Proof of Theorem 5}
 The proof decomposes into three parts. In the first part, we prove that the points selected by the algorithm are everywhere dense in the joint space of $x$ and $\lambda$. In the second part, we show the posterior mean of the GP model converges to the true objective function. Finally, we show that the solution returned by the algorithm converges to the true optimal objective value.

 \subsection{Density of the Design Points}
 We first show that the EI function value for a point $(x', \lambda')$ will be small if it is close to an existing design point. Denote $(x_j, \lambda_j)$ as the closest existing point in the design set to $(x', \lambda')$. Denote $d[(x_j, \lambda_j), (x', \lambda')]$ as the distance between $(x', \lambda')$ and $(x_j, \lambda_j)$. 
 We utilize the Gaussian kernel function $K$ in the GP model and the covariance function can be represented as: $\text{Cov}[ (x_j, \lambda_j),(x', \lambda')] = \sigma_g^2 K((x_j, \lambda_j),(x', \lambda')) =\sigma^2  \exp(-d^2[(x_j, \lambda_j), (x', \lambda')] )$. Here, $d^2[(x, \lambda), (x', \lambda')]=  \sum_{i=1}^{d_x}
 \frac{(x_i-x_i^{\prime})^{2}}{2\theta_{1,i}^2} + \sum_{j=1}^{d_\lambda}\frac{(\lambda_i-\lambda_i^{\prime})^{2}}{2\theta_{2,j}^2}$. $\sigma_g^2$ is the spatial variance of the GP model, $\theta$s are the length scale parameters of the Gaussian kernel, and $d_x$ and $d_\lambda$ are the dimensions of $x$ and $\lambda$, respectively. Then, for any small number $\epsilon>0$, we can find a ball 
 $$B(x_j, \lambda_j|\rho_\epsilon ):=\{(x', \lambda')| d[(x_j, \lambda_j), (x', \lambda')] \leq \rho_\epsilon  \},$$
 around $(x_j, \lambda_j)$ such that all points in $B$ satisfy that
 $$d[(x_j, \lambda_j),(x', \lambda')] \leq \frac{1}{3d_{max}} \min\{\log(1+\epsilon),-\log(1-\epsilon)\},$$
 where $d_{max}$ is the largest value of $d$ between any two points in the joint decision space of $x$ and $\lambda$. For ease of exposition, for any other point $(x, \lambda)$, we denote $d[(x_j, \lambda_j),(x, \lambda)], d[(x_j, \lambda_j),(x', \lambda')], d[(x', \lambda'),(x, \lambda)]$ as $d_1, d_2, d_3$ respectively, and $K[(x_j, \lambda_j),(x, \lambda)], K[(x_j, \lambda_j),(x', \lambda')], K[(x', \lambda'),(x, \lambda)]$ as $K_1, K_2, K_3$ respectively.

 Suppose $(x', \lambda')$ is in $B$. Following the triangular property, we have that
 $$d_3 \leq d_1 + d_2, \ \ d_3\geq d_1 - d_2 .$$
 Thus,
 $$ 
 d_3^2\leq d_1^2+ d_2(2d_1+d_2)\leq d_1^2+3 d_{max} d_2 \leq d_1^2-\log(1-\epsilon),
 $$
 $$ 
 d_3^2\geq d_1^2- d_2(2d_1-d_2) \geq d_1^2-3 d_{max} d_2 \geq d_1^2 -  \log(1+\epsilon).
 $$
 It follows that 
 $$K_3 \geq K_1 (1-\epsilon) \geq  K_1 -\epsilon, \ \ K_3 \leq K_1 (1+\epsilon) \leq  K_1 +\epsilon.$$
 The inequalities hold since $K\leq1$. Hence, $K[(x', \lambda'),(x, \lambda)]  = K[(x_j, \lambda_j),(x, \lambda)] + O(\epsilon)$.

 We now compute the posterior covariance $k_n([(x', \lambda'),(x, \lambda)] )$.
 \begin{equation*}
		 		\begin{aligned}
			 			k_n((x,\lambda),({x'},{\lambda'})) =\sigma_g^2K((x,\lambda),({x'},{\lambda'}))
			 			-c(x,\lambda)^T[R_z+R_\epsilon]^{-1}c({x'},{\lambda'})
			 		\end{aligned}
	 	\label{cov} 
	 \end{equation*}
 Recall that ${R_\epsilon} = \text{Diag}\{\frac{1}{m}\sigma_\epsilon^2\,\cdots,\frac{1}{m}\sigma_\epsilon^2\} =\frac{\sigma_\epsilon^2}{m}  {I}$. When $m \rightarrow \infty$, we have,
 $$[R_z+R_\epsilon]^{-1} = R_z^{-1} + O(\frac{1}{m}) .$$
 Therefore,
 \begin{equation*}
		 		\begin{aligned}
			 			k_n((x,\lambda),({x'},{\lambda'})) &=\sigma_g^2K[(x, \lambda),(x_j, \lambda_j)]
			 			-c(x,\lambda)^TR_z^{-1}c(x_j, \lambda_j) + O(\epsilon) + O(1/m) \\ &= O(\epsilon) + O(1/m). 
			 		\end{aligned}
	 \end{equation*}
 The second equality holds as $(x_j, \lambda_j)$ is a design point and $c(x_j, \lambda_j)$ is the $j$-th column of $R_z$.
 Since $m\rightarrow\infty$, we can select large enough $m$ such that for point $({x'},{\lambda'})$ within $B$, $k_n((x,\lambda),({x'},{\lambda'})) = O(\epsilon)$, $\forall (x,\lambda)$.

 We check the magnitude of $\tilde{\sigma}^2(x'|x', \lambda')$ now. According to Equation (12) in Section 6.3, we can see that $\tilde{\sigma}^2(x'|x', \lambda') = O(\epsilon^{1/2})$. This means the predictive variance can be arbitrary small if $({x'},{\lambda'})$ is close to any existing design point. According to \cite{locatelli1997bayesian} and \cite{wang2022multilevel}, the EI function will tends to 0 when $\tilde{\sigma}^2(x'|x', \lambda') \rightarrow 0$, diminishing the chance of selecting $(x', \lambda')$ for evaluation. 

 On the other hand, if there is no single design points with a ball $B(x_0, \lambda_0|\rho)$, there will be a lower bound for $k_n((x,\lambda),({x_0},{\lambda_0})) $ for points $(x,\lambda)$ in $B$. This lower bound can be computed as if all points outside $B$ were observed accurately. It follows that we can compute a lower bound for $\tilde{\sigma}^2(x_0|x_0, \lambda_0)$ when $N_{mc} \rightarrow \infty$ and thus a lower bound for $\text{EI}(x_0, \lambda_0)$. Therefore, follow the same reasoning as that is \cite{locatelli1997bayesian} and \cite{wang2022multilevel}, the selection rule EI will not keep selecting points around the existing points and remaining certain region unexplored. Instead, it will select points that is everywhere dense in the design space, i.e., for and value $\epsilon>0$, there exist a number $N_\epsilon$, such that for any point $(x_0, \lambda_0)$, there will existing at least one point in $B(x_0, \lambda_0|\epsilon)$ when $N>N_\epsilon$.

 \subsection{Convergence of the Posterior Mean}
 We next show that the posterior mean converges to the true objective function. Following the similar reasoning as \cite{wang2020nonparametric}, with $n$ evaluation points $\{(x_1,\lambda_1),...,(x_n,\lambda_n)\}$ in stage $t$, we can decompose their difference as: 
 \begin{equation*}
		 \left|\mu_n(x)-\tilde{g}_t(x)\right|\leq \left|\mu_n(x)-\mathbb{E}_{\lambda }[Z_n(x,\lambda)]\right|+\left|\mathbb{E}_{\lambda }[Z_n(x,\lambda)]-\tilde{g}_t(x)   \right|.
	 \end{equation*}
 Here, $Z_n$ is the GP posterior for $y(x,\lambda)$ and 
 $$
 \mathbb{E}_{\lambda }[Z_n(x,\lambda)] = \mathbb{E}_{P(\vartheta|\xi^t)} \sum_{l=1}^{\tilde{R}}w_l(\vartheta ) Z_{n}(x, {\lambda}_l(\vartheta)).
 $$

 From the triangle inequality,
 $$
 \mathbb{P}\left\lbrace\left|\mu_n(x)-\tilde{g}_t(x)\right| \geq 2\epsilon \right\rbrace 
 \leq \mathbb{P}\left\lbrace \left|\mu_n(x)-\mathbb{E}_{\lambda}Z_n(x,\lambda)\right|\geq \epsilon \right\rbrace+\mathbb{P}\left\lbrace \left|\mathbb{E}_{\lambda}Z_n(x,\lambda)-\tilde{g}_t(x) \right|\geq \epsilon \right\rbrace.
 $$

 Thus, to show $\mu_n(x)-\tilde{g}_t(x) \to_p 0$, it suffices to establish
 \begin{equation*}
		 \mu_n(x)-\mathbb{E}_{\lambda}Z_n(x,\lambda)\to_p 0,\quad \mbox{and} \quad \mathbb{E}_{\lambda}Z_n(x,\lambda)-\tilde{g}_t(x) \to_p 0
	 \end{equation*}

 We first show the second convergence. Suppose we could directly observe $z(x,\lambda)$ and construct a deterministic Kriging model with conditional mean $\tilde{m}_n(x,\lambda)$ and conditional variance $\tilde{k}_n\left( (x,\lambda),(x,\lambda) \right)$. 
 We have
 $$
 	\left|Z_n(x,\lambda)-z(x,\lambda)\right| \leq \left|Z_n(x,\lambda)-m_n(x,\lambda)  \right|+\left| m_n(x,\lambda)-\tilde{m}_n(x,\lambda)\right|+\left| \tilde{m}_n(x,\lambda)-z(x,\lambda) \right|
 $$
 According to Lemma 5 of \cite{pedrielli2020extended},
 \begin{equation*}
		 \tilde{m}_n(x,\lambda)-m_n(x,\lambda)\to_p 0 \quad \mbox{and} \quad \tilde{k}_n\left( (x,\lambda),(x,\lambda) \right)-k_n\left( (x,\lambda),(x,\lambda) \right)\to_p 0,
	 \end{equation*}
 uniformly as $m$ tends to infinity. Moreover,
 according to Propositions 3.3, 3.4 and 3.5 of \cite{AMS18}, when the design points are everywhere dense, we have
 $
 \tilde{m}_n(x,\lambda)-z(x,\lambda)\to_p 0, \quad \mbox{and} \quad k_n\left( (x,\lambda),(x,\lambda) \right)\to_p 0,
 $
 uniformly for all $(x,\lambda)$.
 In addition, the expectation of $Z_n(x,\lambda)$ is $m_n(x,\lambda)$ and the marginal variance $k_n$ tends to zero uniformly. Thus $\left|Z_n(x,\lambda)-m_n(x,\lambda)  \right|\to_p 0 $ from Chebyshev's Inequality:
 \begin{equation*}
	 \mathbb{P}(\left|Z_n(x,\lambda)-m_n(x,\lambda)  \right|\geq \epsilon)\leq \frac{k_n((x,\lambda),(x,\lambda))}{\epsilon^2}.
	 \end{equation*}
 Combining all that we know, $\left|Z_n(x,\lambda)-z(x,\lambda)\right|\to_p 0$ uniformly for all $(x,\lambda)$.
 And so
 $$
 \sup_{x}\left|\mathbb{E}_{\lambda}Z_n(x,\lambda)-\tilde{g}_t(x)\right|
 =\sup_{x}\left|\mathbb{E}_{\lambda}\left(Z_n(x,\lambda)-z(x,\lambda)\right)\right|
 \leq \sup_{x,\lambda}\left|\left(Z_n(x,\lambda)-z(x,\lambda)\right)\right|\to_p 0,
 $$
 i.e., $ \mathbb{E}_{\lambda}Z_n(x,\lambda)-\tilde{g}_t(x) \to_p 0$ uniformly in $x$. 

 For the first convergence, notice that
 $$
 \left| \mu_n(x)-\mathbb{E}_{\lambda}Z_n(x,\lambda) \right| \leq \left| \mu_n(x)-\mathbb{E}_{\lambda}m_n(x,\lambda)\right|+\left|\mathbb{E}_{\lambda}m_n(x,\lambda)-\mathbb{E}_{\lambda}Z_n(x,\lambda) \right|
 $$
 Since $\mu_n(x)=\frac{1}{N_{mc}}\sum_{i=1}^{N_{mc}} \sum_{l=1}^{\tilde{R}} w_l{(\vartheta_i)} m_n(x, \lambda_l{(\vartheta_i)})$, $\mathbb{E}_{\lambda }[m_n(x,\lambda)] = \mathbb{E}_{\vartheta} \sum_{l=1}^{\tilde{R}}w_l(\vartheta ) m_{n}(x, {\lambda}_l(\vartheta))$, Chebyshev's Inequality implies
 \begin{equation*}
		 \mathbb{P}\left\lbrace \left| \mu_n(x)-\mathbb{E}_{\lambda}m_n(x,\lambda)\right|\geq \epsilon \right\rbrace\leq \frac{\text{Var}_{\lambda} m_n(x,\lambda)}{\epsilon^2N_{mc}}.
	 \end{equation*}
 And so $ \mu_n(x)-\mathbb{E}_{\lambda}m_n(x,\lambda)\to_p 0$ uniformly in $x$ as $N_{mc} \to \infty$.
 On the other hand, 
 \begin{equation*}
		 \sup_{x} \left|\mathbb{E}_{\lambda}m_n(x,\lambda)-\mathbb{E}_{\lambda}Z_n(x,\lambda)\right| \leq \sup_{x,\lambda} \left|  m_n(x,\lambda)-Z_n(x,\lambda) \right| \to_p 0.
	 \end{equation*}
 So $\mathbb{E}_{\lambda}m_n(x,\lambda)-\mathbb{E}_{\lambda}Z_n(x,\lambda)\to_p 0 $. This concludes the proof of the frist convergence.

 \subsection{Convergence of the Optimal Value}
 As the evaluation points we visited are dense, there will exist a subsequence $x_{n_k},\ k=1,2,...$ such that $x_{n_k}\to_p \tilde{x}^*_t$. Denote the best solution identified by the algorithm at $n_k$-th iteration in stage $t$ as $\widehat{x}_{n_k}^*$. For $n_i < N < n_j$, as $n_i$ and $N$ go to infinity, we have
 \begin{equation*}
	 \begin{aligned}
		 \mu_N(\hat{x}_N^*)\leq \mu_N(x_{n_i}) &\leq \tilde{g}_t(x_{n_i}) + o_p(1) = \tilde{g}_t(\tilde{x}^*_t) + o_p(1)
		 \end{aligned}
	 \label{equ:inequality}
	 \end{equation*}
 The last equality is due to the continuous mapping theorem. Therefore,
 \begin{equation*}
	 \tilde{g}_t(\hat{x}_N^*)  \leq \mu_N(\hat{x}_N^*) + o_p(1) \leq \tilde{g}_t(\tilde{x}^*_t) + o_p(1). 	
	 \end{equation*}

 \section{Proof of Theorem 6}
 From Theorem 2, we can see that there exist $t_0$ such that, for any $t>t_0$,
 \begin{equation*}\label{eq: triangle3}
	 P(|\tilde{g}_t(\tilde{x}^*_t)-\tilde{K}(x^*, \vartheta^c)|>\frac{\epsilon}{2})<\frac{\epsilon}{2}.
	 \end{equation*}
 Then for any fixed $t$, we use Theorem 5 to see that there exists
 $N_0$ such that, for any $N>N_0$,
 \begin{equation*}\label{eq: triangle4}
	 P(|\tilde{g}_t(\hat{x}_N^*)-\tilde{g}_t(\tilde{x}^*_t)|>\frac{\epsilon}{2})<\frac{\epsilon}{2}.
	 \end{equation*} 
 Combining these two, the following inequality is immediate:
 \begin{equation*}
	 			P(|\tilde{g}_t(\hat{x}_N^*)-\tilde{K}(x^*, \vartheta^c)|\geq \epsilon)<P(|\tilde{g}_t(\hat{x}_N^*)-\tilde{g}_t(\tilde{x}^*_t)|>\frac{\epsilon}{2})  + P(|\tilde{g}_t(\tilde{x}^*_t)-\tilde{K}(x^*, \vartheta^c)|>\frac{\epsilon}{2}) <\epsilon.
	 \end{equation*}
%
%
%

\bibliographystyle{informs2014} 
\bibliography{sample} 

\begin{thebibliography}{52}
\providecommand{\natexlab}[1]{#1}
\providecommand{\url}[1]{\texttt{#1}}
\providecommand{\urlprefix}{URL }

\bibitem[{Barton et~al.(2014)Barton, Nelson, \protect\BIBand{}
  Xie}]{barton2014quantifying}
Barton RR, Nelson BL, Xie W (2014) Quantifying input uncertainty via simulation
  confidence intervals. \emph{INFORMS Journal on Computing} 26(1):74--87.

\bibitem[{Billingsley(2013)}]{billingsley2013convergence}
Billingsley P (2013) \emph{Convergence of probability measures} (John Wiley \&
  Sons).

\bibitem[{{Bridgewater Associates}(2012)}]{Bridgewater2012}
{Bridgewater Associates} (2012) The all weather story. Technical report,
  Bridgewater Associates, LP,
  \urlprefix\url{https://www.bridgewater.com/research-library/the-all-weather-strategy/}.

\bibitem[{Bui et~al.(2017)Bui, Nguyen, \protect\BIBand{}
  Turner}]{bui2017streaming}
Bui TD, Nguyen C, Turner RE (2017) Streaming sparse gaussian process
  approximations. \emph{Advances in Neural Information Processing Systems} 30.

\bibitem[{Cakmak et~al.(2020)Cakmak, Astudillo~Marban, Frazier,
  \protect\BIBand{} Zhou}]{cakmak2020bayesian}
Cakmak S, Astudillo~Marban R, Frazier P, Zhou E (2020) Bayesian optimization of
  risk measures. \emph{Advances in Neural Information Processing Systems}
  33:20130--20141.

\bibitem[{Costa \protect\BIBand{} Kwon(2019)}]{costa2019risk}
Costa G, Kwon RH (2019) Risk parity portfolio optimization under a markov
  regime-switching framework. \emph{Quantitative finance} 19(3):453--471.

\bibitem[{De~Gunst \protect\BIBand{} Shcherbakova(2008)}]{de2008asymptotic}
De~Gunst M, Shcherbakova O (2008) Asymptotic behavior of bayes estimators for
  hidden markov models with application to ion channels. \emph{Mathematical
  Methods of Statistics} 17:342--356.

\bibitem[{DeMiguel et~al.(2009)DeMiguel, Garlappi, \protect\BIBand{}
  Uppal}]{demiguel2009optimal}
DeMiguel V, Garlappi L, Uppal R (2009) Optimal versus naive diversification:
  How inefficient is the 1/n portfolio strategy? \emph{The review of Financial
  studies} 22(5):1915--1953.

\bibitem[{Dooley et~al.(2010)Dooley, Yan, Mohan, \protect\BIBand{}
  Gopalakrishnan}]{dooley2010inventory}
Dooley KJ, Yan T, Mohan S, Gopalakrishnan M (2010) Inventory management and the
  bullwhip effect during the 2007--2009 recession: evidence from the
  manufacturing sector. \emph{Journal of supply chain management} 46(1):12--18.

\bibitem[{Douc et~al.(2020)Douc, Olsson, \protect\BIBand{}
  Roueff}]{douc2020posterior}
Douc R, Olsson J, Roueff F (2020) Posterior consistency for partially observed
  markov models. \emph{Stochastic Processes and their Applications}
  130(2):733--759.

\bibitem[{Durrett(2019)}]{durrett2019probability}
Durrett R (2019) \emph{Probability: theory and examples}, volume~49 (Cambridge
  university press).

\bibitem[{Fan et~al.(2020)Fan, Hong, \protect\BIBand{}
  Zhang}]{fan2020distributionally}
Fan W, Hong LJ, Zhang X (2020) Distributionally robust selection of the best.
  \emph{Management Science} 66(1):190--208.

\bibitem[{Fox et~al.(2011)Fox, Sudderth, Jordan, \protect\BIBand{}
  Willsky}]{fox2011sticky}
Fox EB, Sudderth EB, Jordan MI, Willsky AS (2011) A sticky hdp-hmm with
  application to speaker diarization. \emph{The Annals of Applied Statistics}
  1020--1056.

\bibitem[{Fr{\"u}hwirth-Schnatter(2006)}]{fruhwirth2006finite}
Fr{\"u}hwirth-Schnatter S (2006) \emph{Finite mixture and Markov switching
  models} (Springer).

\bibitem[{Fu \protect\BIBand{} Healy(1997)}]{fu1997techniques}
Fu MC, Healy KJ (1997) Techniques for optimization via simulation: an
  experimental study on an (s, s) inventory system. \emph{IIE transactions}
  29(3):191--199.

\bibitem[{Hoffman et~al.(2014)Hoffman, Gelman et~al.}]{hoffman2014no}
Hoffman MD, Gelman A, et~al. (2014) The no-u-turn sampler: adaptively setting
  path lengths in hamiltonian monte carlo. \emph{J. Mach. Learn. Res.}
  15(1):1593--1623.

\bibitem[{Hong et~al.(2014)Hong, Hu, \protect\BIBand{} Liu}]{hong2014monte}
Hong LJ, Hu Z, Liu G (2014) Monte carlo methods for value-at-risk and
  conditional value-at-risk: a review. \emph{ACM Transactions on Modeling and
  Computer Simulation (TOMACS)} 24(4):1--37.

\bibitem[{Hong \protect\BIBand{} Nelson(2009)}]{hong2009brief}
Hong LJ, Nelson BL (2009) A brief introduction to optimization via simulation.
  \emph{Proceedings of the 2009 Winter simulation conference (WSC)}, 75--85
  (IEEE).

\bibitem[{Iquebal \protect\BIBand{} Bukkapatnam(2018)}]{iquebal2018change}
Iquebal AS, Bukkapatnam S (2018) Change detection and prognostics for transient
  real-world processes using streaming data. \emph{Recent Advances in
  Optimization and Modeling of Contemporary Problems}, 279--315 (INFORMS).

\bibitem[{Jalali et~al.(2017)Jalali, Van~Nieuwenhuyse, \protect\BIBand{}
  Picheny}]{jalali2017comparison}
Jalali H, Van~Nieuwenhuyse I, Picheny V (2017) Comparison of kriging-based
  algorithms for simulation optimization with heterogeneous noise.
  \emph{European Journal of Operational Research} 261(1):279--301.

\bibitem[{Jones et~al.(1998)Jones, Schonlau, \protect\BIBand{}
  Welch}]{jones1998efficient}
Jones DR, Schonlau M, Welch WJ (1998) Efficient global optimization of
  expensive black-box functions. \emph{Journal of Global optimization}
  13:455--492.

\bibitem[{Keskin et~al.(2022)Keskin, Li, \protect\BIBand{}
  Song}]{keskin2022data}
Keskin NB, Li Y, Song JS (2022) Data-driven dynamic pricing and ordering with
  perishable inventory in a changing environment. \emph{Management Science}
  68(3):1938--1958.

\bibitem[{Lam(2016)}]{lam2016advanced}
Lam H (2016) Advanced tutorial: Input uncertainty and robust analysis in
  stochastic simulation. \emph{2016 Winter Simulation Conference (WSC)},
  178--192 (IEEE).

\bibitem[{Liu \protect\BIBand{} Zhang(2022)}]{liu2022financial}
Liu H, Zhang Y (2022) Financial uncertainty with ambiguity and learning.
  \emph{Management Science} 68(3):2120--2140.

\bibitem[{Liu et~al.(2024)Liu, Lin, \protect\BIBand{} Zhou}]{liu2024bayesian}
Liu T, Lin Y, Zhou E (2024) Bayesian stochastic gradient descent for stochastic
  optimization with streaming input data. \emph{SIAM Journal on Optimization}
  34(1):389--418.

\bibitem[{Locatelli(1997)}]{locatelli1997bayesian}
Locatelli M (1997) Bayesian algorithms for one-dimensional global optimization.
  \emph{Journal of Global Optimization} 10(1):57--76.

\bibitem[{Mankiw(2019)}]{mankiw2019macroeconomics}
Mankiw NG (2019) \emph{Macroeconomics} (New York: Worth Publishers), 10
  edition.

\bibitem[{Mor et~al.(2021)Mor, Garhwal, \protect\BIBand{}
  Kumar}]{mor2021systematic}
Mor B, Garhwal S, Kumar A (2021) A systematic review of hidden markov models
  and their applications. \emph{Archives of computational methods in
  engineering} 28:1429--1448.

\bibitem[{Oprisor \protect\BIBand{} Kwon(2020)}]{oprisor2020multi}
Oprisor R, Kwon R (2020) Multi-period portfolio optimization with investor
  views under regime switching. \emph{Journal of Risk and Financial Management}
  14(1):3.

\bibitem[{Park \protect\BIBand{} Gupta(2011)}]{park2011regime}
Park S, Gupta S (2011) A regime-switching model of cyclical category buying.
  \emph{Marketing Science} 30(3):469--480.

\bibitem[{Pearce \protect\BIBand{} Branke(2017)}]{pearce2017bayesian}
Pearce M, Branke J (2017) Bayesian simulation optimization with input
  uncertainty. \emph{2017 Winter Simulation Conference (WSC)}, 2268--2278
  (IEEE).

\bibitem[{Pedrielli et~al.(2020)Pedrielli, Wang, \protect\BIBand{}
  Ng}]{pedrielli2020extended}
Pedrielli G, Wang S, Ng SH (2020) An extended two-stage sequential optimization
  approach: Properties and performance. \emph{European Journal of Operational
  Research} .

\bibitem[{Pun et~al.(2023)Pun, Wang, \protect\BIBand{} Yan}]{pun2023data}
Pun CS, Wang T, Yan Z (2023) Data-driven distributionally robust cvar portfolio
  optimization under a regime-switching ambiguity set. \emph{Manufacturing \&
  Service Operations Management} 25(5):1779--1795.

\bibitem[{Royden \protect\BIBand{} Fitzpatrick(2010)}]{royden2010real}
Royden H, Fitzpatrick PM (2010) \emph{Real analysis} (China Machine Press).

\bibitem[{Rudin et~al.(1964)}]{rudin1964principles}
Rudin W, et~al. (1964) \emph{Principles of mathematical analysis}, volume~3
  (McGraw-hill New York).

\bibitem[{Song \protect\BIBand{} Shanbhag(2019)}]{song2019stochastic}
Song E, Shanbhag UV (2019) Stochastic approximation for simulation optimization
  under input uncertainty with streaming data. \emph{2019 Winter Simulation
  Conference (WSC)}, 3597--3608 (IEEE).

\bibitem[{Stuart \protect\BIBand{} Teckentrup(2018)}]{AMS18}
Stuart AM, Teckentrup AL (2018) {Posterior consistency for Gaussian process
  approximations of Bayesian posterior distributions}. \emph{Mathematics of
  Computation of the American Mathematical Society} 87(310):721--753.

\bibitem[{Trapero et~al.(2015)Trapero, Kourentzes, \protect\BIBand{}
  Fildes}]{trapero2015identification}
Trapero JR, Kourentzes N, Fildes R (2015) On the identification of sales
  forecasting models in the presence of promotions. \emph{Journal of the
  operational Research Society} 66(2):299--307.

\bibitem[{Wang et~al.(2020{\natexlab{a}})Wang, Yuan, \protect\BIBand{}
  Ng}]{wang2020gaussian}
Wang H, Yuan J, Ng SH (2020{\natexlab{a}}) Gaussian process based optimization
  algorithms with input uncertainty. \emph{IISE Transactions} 52(4):377--393.

\bibitem[{Wang et~al.(2020{\natexlab{b}})Wang, Zhang, \protect\BIBand{}
  Ng}]{wang2020nonparametric}
Wang H, Zhang X, Ng SH (2020{\natexlab{b}}) A nonparametric bayesian approach
  for simulation optimization with input uncertainty. \emph{arXiv preprint
  arXiv:2008.02154} .

\bibitem[{Wang et~al.(2022)Wang, Ng, \protect\BIBand{}
  Haskell}]{wang2022multilevel}
Wang S, Ng SH, Haskell WB (2022) A multilevel simulation optimization approach
  for quantile functions. \emph{INFORMS Journal on Computing} 34(1):569--585.

\bibitem[{Wang et~al.(2024)Wang, Wang, Xia, \protect\BIBand{}
  Xu}]{wang2024simulation}
Wang S, Wang H, Xia J, Xu X (2024) Simulation optimization with non-stationary
  streaming input data. \emph{Proceedings of the Winter Simulation Conference},
  3494--3505.

\bibitem[{Wang(2025)}]{wang2025limited}
Wang Y (2025) Limited firm insurance and aggregate implications.
  \emph{Management Science} 71(7):5997--6046.

\bibitem[{Wang \protect\BIBand{} Zhou(2025)}]{wang2025ranking}
Wang Y, Zhou E (2025) Ranking and selection with simultaneous input data
  collection. \emph{arXiv preprint: 2503.11773} .

\bibitem[{Wilson(2007)}]{wilson2007impact}
Wilson MC (2007) The impact of transportation disruptions on supply chain
  performance. \emph{Transportation Research Part E: Logistics and
  Transportation Review} 43(4):295--320.

\bibitem[{Wu et~al.(2024)Wu, Wang, \protect\BIBand{} Zhou}]{wu2024data}
Wu D, Wang Y, Zhou E (2024) Data-driven ranking and selection under input
  uncertainty. \emph{Operations Research} 72(2):781--795.

\bibitem[{Wu \protect\BIBand{} Zhou(2017)}]{wu2017ranking}
Wu D, Zhou E (2017) Ranking and selection under input uncertainty: A budget
  allocation formulation. \emph{2017 Winter Simulation Conference (WSC)},
  2245--2256 (IEEE).

\bibitem[{Wu et~al.(2018)Wu, Zhu, \protect\BIBand{} Zhou}]{wu2018bayesian}
Wu D, Zhu H, Zhou E (2018) A bayesian risk approach to data-driven stochastic
  optimization: Formulations and asymptotics. \emph{SIAM Journal on
  Optimization} 28(2):1588--1612.

\bibitem[{Xiao \protect\BIBand{} Gao(2018)}]{xiao2018simulation}
Xiao H, Gao S (2018) Simulation budget allocation for selecting the top-m
  designs with input uncertainty. \emph{IEEE Transactions on Automatic Control}
  63(9):3127--3134.

\bibitem[{Xie et~al.(2014)Xie, Nelson, \protect\BIBand{}
  Barton}]{xie2014bayesian}
Xie W, Nelson BL, Barton RR (2014) A bayesian framework for quantifying
  uncertainty in stochastic simulation. \emph{Operations Research}
  62(6):1439--1452.

\bibitem[{Yee Whye~Teh \protect\BIBand{} Blei(2006)}]{Teh01122006}
Yee Whye~Teh MJB Michael I~Jordan, Blei DM (2006) Hierarchical dirichlet
  processes. \emph{Journal of the American Statistical Association}
  101(476):1566--1581.

\bibitem[{Zhou \protect\BIBand{} Xie(2015)}]{zhou2015simulation}
Zhou E, Xie W (2015) Simulation optimization when facing input uncertainty.
  \emph{2015 Winter Simulation Conference (WSC)}, 3714--3724 (IEEE).

\end{thebibliography}


\end{document}